\documentclass[11pt, leqno]{amsart}
\usepackage{amsthm,amsmath,amssymb,amscd,verbatim}
\usepackage[curve,matrix,arrow,frame,tips]{xy}
\theoremstyle{plain}
\newtheorem{theorem}{Theorem}[section]
\newtheorem{lemma}[theorem]{Lemma}
\newtheorem{cor}[theorem]{Corollary}
\newtheorem{prop}[theorem]{Proposition}

\newtheorem{thm}{Theorem}[section]
\newtheorem{coro}[thm]{Corollary}

\theoremstyle{remark}

\theoremstyle{definition}

\newtheorem{defn}[theorem]{Definition}
\newtheorem{Remark}[theorem]{Remark}

\numberwithin{equation}{subsection}
\numberwithin{theorem}{subsection}

\newcommand\hfld[2]{\smash{\mathop{\hbox to 10mm{\rightarrowfill}}
     \limits^{\scriptstyle#1}_{\scriptstyle#2}}}
\newcommand\hflg[2]{\smash{\mathop{\hbox to 10mm{\leftarrowfill}}
     \limits^{\scriptstyle#1}_{\scriptstyle#2}}}

\title{Base change for Bernstein centers of \\ depth zero principal series blocks}
\author{Thomas J. Haines}

\date{}
\begin{document}

\thanks{Research partially supported by NSF grants FRG-0554254, DMS-0901723, and a University of Maryland GRB Semester Award.}

\date{}

\maketitle

\begin{abstract} \begin{small}
Let $G$ be an unramified group over a $p$-adic field.  This article introduces a base change homomorphism for Bernstein centers of depth-zero principal series blocks for $G$ and proves the corresponding base change fundamental lemma.  This result is used in the approach to Shimura varieties with $\Gamma_1(p)$-level structure initiated by M.~Rapoport and the author in \cite{HR2}. 
\end{small} 
\end{abstract} 
\markboth{T. Haines}
{Base change for depth zero principal series blocks}


\section{Introduction}

Let $F$ denote a $p$-adic field, and let $F_r \supset F$ denote the unique degree $r$ unramified extension of $F$ contained in some algebraic closure $\bar{F}$ of $F$.  Let $\theta$ denote a generator of ${\rm Gal}(F_r/F)$.  Let $G$ denote an unramified connected reductive group over $F$.  The automorphism $\theta$ determines an automorphism of $G(F_r)$, which we also denote with the symbol $\theta$.

Using $\theta$, we have the notion of stable twisted orbital integral ${\rm SO}_{\delta \theta}(\phi)$ for any locally constant compactly-supported function $\phi$ on $G(F_r)$ and any element $\delta \in G(F_r)$  with semisimple norm.  See \cite{Ko82} for the definition of the norm map $\mathcal N$ from stable $\theta$-conjugacy classes in $G(F_r)$ to stable conjugacy classes in $G(F)$.  For a precise definition of ${\rm SO}_{\delta \theta}$, see e.g.~\cite{Ko86b}, or \cite{H09a}, (5.1.2).

This article is concerned with the matching of the (twisted) orbital integrals of certain 
functions on the groups $G(F_r)$ and $G(F)$, respectively.  If $\phi \in \mathcal H(G(F_r))$ and $f \in \mathcal H(G(F))$ are functions in the corresponding Hecke algebras of locally constant compactly-supported functions, then we say $\phi, f$ are {\em associated} (or {\em have matching orbital integrals}), if the following result holds for the stable (twisted) orbital integrals:  for every semisimple element $\gamma \in G(F)$, we have
$$
{\rm SO}_\gamma(f) = \sum_{\rm \delta}\Delta(\gamma, \delta) {\rm SO}_{\delta \theta}(\phi)
$$
where the sum is over stable $\theta$-conjugacy classes $\delta \in G(F_r)$ with semisimple norm, and where $\Delta(\gamma, \delta) = 1$ if $\mathcal N\delta = \gamma$ and $\Delta(\gamma,\delta) = 0$ otherwise.  See e.g. \cite{Ko86b}, \cite{Ko88}, \cite{Cl90}, or \cite{H09a} for further discussion.  

Of primary importance is the case of spherical Hecke algebras.  Suppose $K_r \subset G(F_r)$ and $K \subset G(F)$ are hyperspecial maximal compact subgroups associated to a hyperspecial vertex in the Bruhat-Tits building $\mathcal B(G(F))$ for $G(F)$, and suppose $\phi \in \mathcal H(G(F_r))$ belongs to the corresponding spherical Hecke algebra $\mathcal H_{K_r}(G(F_r))$.  The Satake isomorphism gives rise to a natural algebra homomorphism
$$
b_r: \mathcal H_K(G(F_r)) \rightarrow \mathcal H_K(G(F)),
$$
cf.~\cite{Cl90}.  The base change fundamental lemma for spherical functions asserts that $\phi$ and $b_r(\phi)$ are associated.  This was proved by Clozel \cite{Cl90} and Labesse \cite{Lab90}.  Even the earlier special cases of ${\rm GL}_2$ \cite{Lan} and ${\rm GL}_n$ \cite{AC} gave rise to important global and local applications, such as the existence of base-change lifts of certain automorphic representations for ${\rm GL}_n$.  The theorem of Clozel and Labesse played an important role in Kottwitz' work \cite{Ko90, Ko92, Ko92b} on Shimura varieties with good reduction at $p$.

In \cite{H09a} an analogous base change fundamental lemma is proved for centers of parahoric Hecke algebras.  It plays a role in the study of Shimura varieties with parahoric level structure at $p$, 
see \cite{H05}.  
A very special case of \cite{H09a} relates to the center of the Iwahori-Hecke algebra for $G(F)$.   By Bernstein's theory \cite{BD}, this can be viewed as the ring of regular functions on the variety of supercuspidal supports of the Iwahori block (the subcategory $\mathcal R_{I}(G)$ of the category $\mathcal R(G)$ of smooth representations of $G(F)$ whose objects are generated by their invariants under an Iwahori subgroup).

The purpose of this article is to generalize this result to certain other principal series blocks in the Bernstein decomposition, namely the {\em depth-zero} principal series blocks.   

To state the theorem, we need some more notation.  Denote the ring of integers of $F$ by $\mathcal O$ and the residue field by $k$.   Let $A$ denote a maximal $F$-split torus, and set $T := {\rm Cent}_G(A)$, a maximal torus of $G$ defined and unramified over $F$.  Now choose an Iwahori subgroup $I \subset G(F)$ which is in good position relative to $T$, that is, the alcove ${\bf a}$ in $\mathcal B(G(F))$ corresponding to $I$ is contained in the apartment $\mathcal A_T$ corresponding to $A$.  Let $I^+$ denote the pro-unipotent radical of $I$.

Let $T(F)_1$ denote the maximal compact open subgroup of $T(F)$, and let $T(F)_1^+ = T(F)_1 \cap I^+$ denote its pro-unipotent radical.  Throughout this article, $\chi$ will denote a depth-zero character on $T(F)_1$.  This means that $\chi$ factors through a character $T(F)_1/T(F)_1^+ \rightarrow \mathbb C^\times$, which we also denote by $\chi$.  Via the canonical isomorphism
$$
T(F)_1/T(F)_1^+  ~ \widetilde{\rightarrow} ~ I/I^+,
$$
we see that $\chi$ induces a smooth character $\rho := \rho_\chi$ on $I$ which is trivial on $I^+$.   Then we may consider the Hecke algebra $\mathcal H(G,\rho) := \mathcal H(G,I,\chi)$, which is defined as   
$$
\mathcal H(G,\rho) = \{ f \in \mathcal H(G) ~ | ~ f(i_1 g i_2) = \rho^{-1}(i_1) f(g) \rho^{-1}(i_2) \,\,\ \forall i_1,i_2 \in I, \,\forall g \in G \}.
$$
(Cf. \cite{Ro}.)  Convolution is defined using the Haar measure which gives $I$ volume 1.  Write $\mathcal Z(G,\rho)$ for the center of $\mathcal H(G,\rho)$.  

Let $\frak s = \frak s_\chi$ be the inertial equivalence class of the cuspidal pair $(T(F),\widetilde{\chi})_G$, for any extension of $\chi$ to a character $\widetilde{\chi}: T(F) \rightarrow \mathbb C^\times$.  Then $\frak s$ depends only on the relative Weyl group orbit of $\chi$.  We are concerned with the principal series Bernstein block
$$
\mathcal R_{\mathfrak s}(G) = \mathcal R_\chi(G),
$$
the category of smooth representations of $G(F)$ whose irreducible objects are constituents of some normalized principal series $i^G_B(\xi)$, where $\xi:T(F) \rightarrow \mathbb C^\times$ has $\xi|_{T(F)_1} = \chi$.  It turns out (cf. Proposition \ref{type}) that $(I, \rho_\chi)$ is a Bushnell-Kutzko type for 
$\mathcal R_\chi(G)$, so that the Bernstein center of $\mathcal R_\chi(G)$ can be identified with the ring $\mathcal Z(G,\rho)$.

Let $N_r:T(F_r)_1 \rightarrow T(F)_1$ denote the norm homomorphism given by $t \mapsto t\, \theta(t) \cdots 
\theta^{r-1}(t)$.  The character $\chi_r:= \chi \circ N_r$ is a depth-zero character on $T(F_r)_1$, and gives rise to the Bernstein block $\mathcal R_{\chi_r}(G_r)$ for the $p$-adic group $G_r$, the character $\rho_r$ on $I_r$, and the Hecke algebra $\mathcal H(G_r,\rho_r)$ with center $\mathcal Z(G_r,\rho_r)$.  Here $G_r := G(F_r)$ and $I_r \subset G_r$ is the Iwahori subgroup corresponding to $I$.

In Definition \ref{bc_defn} we define a {\em base change homomorphism} 
\begin{equation} \label{bc_chi_defn}
b_r: \mathcal Z(G_r, \rho_r) \rightarrow \mathcal Z(G,\rho).
\end{equation}
This is analogous to the base change homomorphisms defined for spherical Hecke algebras or centers of parahoric Hecke algebra (cf. \cite{H09a}).  Our main theorem is the following.

\begin{thm} \label{main_thm}
If $\phi \in \mathcal Z(G_r,\rho_r)$, the functions $\phi$ and $b_r(\phi)$ are associated.
\end{thm}

Now write $\mathcal H(G,I^+)$ for $\mathcal H_{I^+}(G(F))$ and $\mathcal Z(G,I^+)$ for its center.  In section \ref{bc_I+_sec}, we use (\ref{bc_chi_defn}) to define a natural algebra homomorphism 
\begin{equation} \label{bc_I+_defn0}
b_r: \mathcal Z(G_r, I^+_r) \rightarrow \mathcal Z(G,I^+).
\end{equation}
Moreover, we show how Theorem \ref{main_thm} immediately implies the following result.

\begin{coro} \label{bc_I+_cor}
If $\phi \in \mathcal Z(G_r,I^+_r)$, the functions $\phi$ and $b_r(\phi)$ are associated.
\end{coro}

In \cite{HR2}, these results are used in the special case of $G = {\rm GL}_d$ to study Shimura varieties in the Drinfeld case with $\Gamma_1(p)$-level structure at $p$.  In future works they will be applied to other Shimura varieties with $\Gamma_1(p)$-level structure.  Here, by ``$\Gamma_1(p)$-level'' we mean that the compact open subgroup $K_p$ coming from the Shimura data $({\bf G},X,K_pK^p)$ is the pro-unipotent radical of an Iwahori subgroup of ${\mathbf G}(\mathbb Q_p)$.  Theorem \ref{main_thm} and Corollary \ref{bc_I+_cor} will play a role in the pseudostabilization of the counting points formula for the Shimura varieties just mentioned (see \cite{HR2} for the Drinfeld case).

\medskip

This article overlaps with \cite{H09a} in the case where $\chi = {\rm triv}$, for then $\mathcal H(G,\rho)$ is just the Iwahori-Hecke algebra and the theorem proved here is a special case of \cite{H09a}.  Here we use Labesse's method of elementary functions \cite{Lab90}, whereas in \cite{H09a} we followed Clozel's method \cite{Cl90} more closely.  When $\chi = {\rm triv}$, many arguments presented here become simpler.  Thus in that special case, this article gives an alternative (somewhat easier) proof for the Iwahori-Hecke algebra special case of \cite{H09a}.

As in earlier papers, the proof of Theorem \ref{main_thm} is by induction on the semisimple rank of $G$, and so consists of two steps: (i) use descent formulas to reduce to the case of elliptic elements, and (ii) use a global trace formula argument to prove the theorem for suitable elliptic elements.  Labesse's elementary functions are used in step (ii) to give sufficiently many character identities that the required character identity for central elements is forced (see section \ref{ell_proof}). 

Here is an outline of the contents of the paper.  In section \ref{notation_subsec} we recall some standard notation which will be used throughout the article.  In section \ref{Bern_sec} we review the essential facts about the Bernstein decomposition related to depth-zero principal series blocks.  In section \ref{BC_sec} we define the base change homomorphism and prove some of its basic properties.  
The constant term homomorphism (an essential ingredient for descent) is studied in section \ref{CT_sec}, and in particular is proved to be compatible with the base change homomorphisms.  The descent formulas themselves are the subject of section \ref{Des_sec}.  Section \ref{reductions_section} reduces the general fundamental lemma to the case where $G$ is adjoint and $\gamma$ is a norm and is elliptic and strongly regular semisimple.  (Unfortunately, this section is the most technical section of all.)  
In section \ref{elem_fcns_sec} we introduce and study suitable analogues of Labesse's elementary functions adapted to the Bernstein component $\mathcal R_\chi(G)$.  In section \ref{ell_proof} we use all the preceding material to conclude the proof by establishing the character identity which is equivalent, by the existence of the local data, to the required identity of stable (twisted) orbital integrals.  In section \ref{bc_I+_sec} we define (\ref{bc_I+_defn0}) and explain how Theorem \ref{main_thm} implies Corollary \ref{bc_I+_cor}.  Finally, in section \ref{corr_sec} we correct and clarify a few minor mistakes in \cite{H09a}.

\medskip

\noindent {\em Acknowledgments}: I thank Alan Roche for helpful conversations.  I thank Michael Rapoport for his interest in this work.  I thank especially Robert Kottwitz for providing crucial help with the proof of Lemma \ref{C_T_lem}.  Much of this work was written during Fall 2010 at the Institute for Advanced Study in Princeton.  I thank the IAS for providing financial support\footnote{This material is based upon work supported by the National Science Foundation under agreement No. DMS-0635607.  Any opinions, findings and conclusions or recommendations expressed in this material are those of the author and do not necessarily reflect the views of the National Science Foundation.} and an excellent working environment.

\section{Further notation} \label{notation_subsec}

Let $L$ denote $\widehat{F^{un}}$, the completion of the maximal unramified extension of $F$ contained in $\bar{F}$.  Let $F_r/F$ be an unramified extension of degree $r$ contained in $L$, with ring of integers $\mathcal O_r$ and residue field $k_r$.  Let $\sigma \in {\rm Aut}(L/F)$ denote the Frobenius automorphism, and use the same symbol to denote the induced automorphism on groups of the form $G(L)$, etc.  Fix an algebraic closure $\bar{L}$ for $L$ and define the inertia subgroup\footnote{Not to be confused with the Iwahori subgroup!} as $I = {\rm Gal}(\bar{L}/L)$.

Denote groups of $F_r$-points with a subscript $r$, e.g. $T_r := T(F_r)$ and $G_r := G(F_r)$.  Fix a generator $\theta \in {\rm Gal}(F_r/F)$.  We use the same symbol $\theta$ to denote the induced automorphisms of groups of $F_r$-points $T_r$, $G_r$, etc.

Let $\varpi$ denote a uniformizer for the field $F$.  

We recall the basic facts on the Kottwitz homomorphism \cite{Ko97}.  In loc.~cit. is defined a canonical surjective homomorphism for any connective reductive $F$-group $H$
\begin{equation} \label{K-hom}
\kappa_H: H(L) \twoheadrightarrow X^*(Z(\widehat{H}))_I,
\end{equation}
where $\widehat{H} = \widehat{H}(\mathbb C)$ denotes the Langlands dual group of $H$.  
By loc.~cit., it remains surjective on taking $\sigma$-fixed points:
$$
\kappa_H: H(F) \twoheadrightarrow X^*(Z(\widehat{H}))_I^\sigma.
$$
We define
\begin{align*}
H(L)_1 &:= {\rm ker}(\kappa_H) \\
H(F)_1 &:= {\rm ker}(\kappa_H) \cap H(F).
\end{align*}
If $H$ is any {\em unramified} $F$-torus, then $H(F)_1$ is indeed the unique maximal compact subgroup of $H(F)$, as the notation in the introduction indicated (cf.~\cite{H09a}, Lemma 2.4.1). 

Let $N_r$ denote the norm homomorphism $T(F_r)_1 \rightarrow T(F)_1$, given by
$$
N_r(t) = t \, \theta(t) \, \cdots \, \theta^{r-1}(t).$$
We will use the same symbol to denote the norm homomorphisms $N_r: T(F_r) \rightarrow T(F)$ and $N_r: T(F_r)_1/T(F_r)_1^+ \rightarrow T(F)_1/T(F)_1^+$.

Let $\chi_r := \chi \circ N_r$.  This will be thought of as a character on $T(F_r)_1/T(F_r)_1^+$ or as a depth-zero character on $T(F_r)_1$, depending on context.

The set of all smooth characters on $T(F)$ carries a natural (left) action under the relative Weyl group $W(F) = N_GT(F)/T(F)$; let $W_\chi$ denote the subgroup of $W(F)$ which fixes $\chi$.  Likewise we define $W_{\chi_r}$ in the relative Weyl group $W_r = W(F_r)$.  

Fix an $F$-rational Borel subgroup $B$ containing $T$ and a Levi decomposition $B = TU$, where $U$ is the unipotent radical of $B$.  There is an Iwahori decomposition of $I$ with respect to $T$ and any such $B$:
\begin{equation} \label{Iwahori_decomp}
I = I_U \cdot I_T \cdot I_{\bar{U}},
\end{equation}
where $I_U := U \cap I$, etc..  Note that $I_T = T(F)_1$.  

Let $H$ denote any group.  For $h \in H$ and any subset $S \subset H$, we set $^hS := hSh^{-1}$.  For a function $f$ on $H$, we define a new function $^hf$ by $^hf(x) := f(h^{-1}xh)$ for $x \in H$. 

For two smooth $G(F)$-representations $\pi$ and $\pi'$, the notation $\pi \in \pi'$ will mean that $\pi$ is isomorphic to a subquotient of $\pi'$.

\section{Bernstein center for depth zero principal series} \label{Bern_sec}

We assume $G$ is any connected reductive group which is unramified over $F$.

\subsection{Review of depth zero principal series blocks}

When the base field $F$ is understood, we will often abbreviate $G(F)$ by the symbol $G$ and $W(F)$ by the symbol $W$. 

Let $\chi$ denote a depth-zero character $\chi: T(F)_1 \rightarrow \mathbb C^\times$.  Let $\widetilde{\chi}$ denote any extension of $\chi$ to a character $\widetilde{\chi}: T(F) \rightarrow \mathbb C^\times$.  We consider the inertial class 
$$\mathfrak s = \mathfrak s_\chi = [T(F),\widetilde{\chi}]_G.
$$
The inertial class $\mathfrak s$ depends only on the $W$-orbit of $\chi$.  

Let $\mathcal R(G)$ denote the category of smooth representations of $G$.  Let $\mathcal R_\mathfrak s(G)$ denote the Bernstein component indexed by $\mathfrak s$.  Recall that this is the full subcategory of $\mathcal R(G)$ whose objects have the property that each of their subquotients is a subquotient of a principal series representation $i^G_B(\widetilde{\chi}\eta)$, for some unramified character $\eta$ of $T(F)$.  We shall often denote $\mathcal R_{\mathfrak s}(G)$ by $\mathcal R_\chi(G)$.

\subsection{Algebraic variety associated to a depth-zero block}

Fix $\chi$ and ${\mathfrak s} = {\mathfrak s}_\chi$ as above, and the associated Bernstein component $\mathcal R_{\chi}(G)$.   Recall (cf. e.g.~\cite{HR2}, $\S9.2$) the set 
$$
\mathfrak X_{\mathfrak s} =  \{ (T,\xi)_G \}
$$ 
of {\em supercuspidal supports} $(T,\xi)_G$ of irreducible representations $\pi$ in the category $\mathcal R_{\chi}(G)$.  Here $\xi:T(F)\rightarrow \mathbb C^\times$ is a smooth character extending some $W(F)$-conjugate of $\chi$ and $(T,\xi)_G$ denotes the $G$-conjugacy class of the pair $(T,\xi)$.  Recall that $\chi$ possesses at least one $W_\chi$-invariant extension $\widetilde{\chi}$ (see \cite{HR2}, Remark 9.2.4\footnote{This handles the split case; the unramified case is similar.}), and that having fixed such a $\widetilde{\chi}$, we have a bijection
\begin{align*}
\widehat{A}/W_\chi ~ &\widetilde{\rightarrow} ~ \mathfrak X_{\mathfrak s} \\
\eta &\mapsto (T,\widetilde{\chi} \eta)_G,
\end{align*}
where $\eta \in \widehat{A}$ is viewed as an unramified character on $T(F)$ (cf. e.g. \cite{H09a}, Lemma 2.4.2).    We use this bijection to endow $\mathfrak X_{\mathfrak s}$ with the structure of an affine algebraic variety.  Up to isomorphism, this structure depends neither on the choice of $\chi$ in its $W(F)$-orbit, nor on the choice of the extension $\widetilde{\chi}$ of $\chi$.  We have
$$
\mathfrak X_\mathfrak s = {\rm Spec}(\mathbb C[{\mathfrak X_\mathfrak s}]),
$$
for a $\mathbb C$-algebra $\mathbb C[{\mathfrak X_{\mathfrak s}}]$ which is isomorphic to $\mathbb C[X_*(A)]^{W_\chi}$.  

We will often denote $\mathfrak X_{\mathfrak s}$ by $\mathfrak X_\chi$.

\subsection{The associated type for $\mathcal R_\chi(G)$}

Let $\rho = \rho_\chi: I \rightarrow \mathbb C^\times$ be the smooth character which is induced by $\chi: I/I^+ =  T(F)_1/T(F)_1^+ \rightarrow \mathbb C^\times$. 

\begin{prop} \label{type}
The pair $(I,\rho)$ is a Bushnell-Kutzko type for $\mathcal R_\chi(G)$.
\end{prop}

By definition, this means that an irreducible representation $\pi \in \mathcal R(G)$ belongs to $\mathcal R_{\chi}(G)$ if and only if its restriction to $I$ contains the representation $\rho$.  See \cite{BK}.

\begin{proof} ({\em Sketch}.)  
For $F$-split groups $G$, Roche \cite{Ro} produces a Bushnell-Kutzko type for $\mathcal R_\chi(G)$ for {\em any} smooth character $\chi$, under certain restrictions on the residual characteristic of $F$ (and in the case of depth-zero $\chi$, his type is $(I,\rho)$).  This relies on his theory of Hecke algebra isomorphisms and hence involves, in the case of general $\chi$, some restrictions on the residual characteristic of $F$.

When $\chi$ is depth-zero, but now assuming only that $G$ is unramified over $F$, one can prove more directly that $(I,\rho)$ is a type for $\mathcal R_\chi(G)$.  One proves that for any irreducible representation $\pi \in \mathcal R(G)$, the Jacquet functor $\pi \mapsto \pi_U$ induces a $T(F)_1$-equivariant isomorphism
\begin{equation} \label{Jac_isom}
\pi^\rho ~ \widetilde{\rightarrow} ~ \pi_U^\chi.
\end{equation}
One then deduces $(I,\rho)$ is a type for $\mathcal R_\chi(G)$ using  Frobenius reciprocity.  The proof of (\ref{Jac_isom}) relies on a consequence of certain Hecke algebra isomorphisms, namely: any non-zero element of $\mathcal H(G,\rho)$ supported on a single $I$-double coset is invertible.  But for the depth-zero characters such isomorphisms were established by Goldstein \cite{Go} (in the split case) and by Morris \cite{Mor}, Theorem 7.12 (in general), and these are valid with no restrictions on the residual characteristic.  Details for this approach in the split case are contained in \cite{H09b}, and the unramified case is handled in exactly the same way.
\end{proof}

\subsection{The action of the center $\mathcal Z(G,\rho)$ of $\mathcal H(G,\rho)$}

Proposition \ref{type} means that $\mathcal R_\chi(G)$ is naturally equivalent to the category of $\mathcal H(G,\rho)$-modules (cf. \cite{Ro}).  It also means that there is a canonical algebra isomorphism
\begin{equation} \label{isom1}
\beta: \mathbb C[{\mathfrak X_{\mathfrak s}}] ~ \widetilde{\rightarrow} ~ \mathcal Z(G,\rho)
\end{equation}
which is characterized as follows.  For an extension $\xi$ of some $W$-conjugate of $\chi$, consider the space $i^G_B(\xi)^{\rho}$ of locally constant functions $f: G \rightarrow \mathbb C$ such that
$$
f(tugi) = (\delta_B^{1/2}\xi)(t)f(g)\rho(i) 
$$
for all $t\in T(F)$, $u \in U(F)$, $g \in G$, and $i \in I$.  Then 
$$
z \in \mathcal Z(G,\rho) \,\, \mbox{\em acts on the left on $i^G_B(\xi)^\rho$ by the 
scalar $\beta^{-1}(z)(\xi)$}.
$$  
Here $\beta^{-1}(z)$ is viewed as a regular function on the variety $\mathfrak X_{\mathfrak s}$ and the argument $\xi$ is an abbreviation for the point $(T,\xi)_G \in \mathfrak X_{\mathfrak s}$.  

\begin{Remark} \label{^wbeta} Note that $\beta^{-1}(z)$ being well-defined as a function on the class $(T,\xi)_G$ means that $\beta^{-1}(z)(\, ^w\xi) = \beta^{-1}(z)(\xi)$ for all $w \in W$; that is, $\beta^{-1}(z)$ is $W$-invariant as a function of $\xi$.
\end{Remark}

\subsection{Right vs.~left actions on principal series}

For later use (e.g. in the proof of Proposition \ref{pres_center}), we need to rephrase this in terms of the {\em right} action of $\mathcal Z(G,\rho)$ on principal series representations.  Let $\iota: \mathcal H(G) \rightarrow \mathcal H(G)$ denote the involution of the Hecke algebra of locally constant compactly-supported functions $\mathcal H(G)$, defined by the identity
$$
(\iota h)(g) = h(g^{-1}) \,\,\,, \,\,\,  h \in \mathcal H(G),\,  g \in G.
$$
We use $\iota$ to convert left $\mathcal H(G)$-actions into right $\mathcal H(G)$-actions, and vice-versa.  It is clear that $\mathcal H(G,\rho^{-1})$ (resp. $\mathcal H(G,\rho)$) acts on the left (resp. right) on the space 
$i^G_B(\xi^{-1})^{\rho^{-1}}$, and that for $\Psi \in i^G_B(\xi^{-1})^{\rho^{-1}}$, and $h \in \mathcal H(G,\rho)$, the left and right actions are related by the formula
\begin{equation} \label{iota(h)}
\iota h \cdot \Psi = \Psi \cdot h.
\end{equation}

Duality gives a perfect pairing
$$
(\cdot \, , \cdot): i^G_B(\xi)^{\rho} \times i^G_B(\xi^{-1})^{\rho^{-1}} \rightarrow \mathbb C.
$$
For $z \in \mathcal Z(G,\rho)$, $\phi_1 \in i^G_B(\xi)^{\rho}$, and $\phi_2 \in i^G_B(\xi^{-1})^{\rho^{-1}}$, we have the relation
$$
(z(\phi_1), \phi_2) = (\phi, \iota z(\phi_2)).
$$
Thus
\begin{equation} \label{beta/iota}
\beta^{-1}(z)(\xi) = \beta^{-1}(\iota z)(\xi^{-1}).
\end{equation}
From (\ref{iota(h)}) with $h=z \in \mathcal Z(G,\rho)$ and (\ref{beta/iota}), we deduce the following result.

\begin{lemma} \label{ch_xi_lem}
Let the Haar measure on $G(F)$ be normalized so that $I$ has measure 1.  Then $z \in \mathcal Z(G,\rho)$ acts via right convolutions on $i^G_B(\xi^{-1})^{\rho^{-1}}$ by the scalar
$$
ch_{\xi^{-1}}(z) := \beta^{-1}(z)(\xi).
$$
\end{lemma}

\section{Base change homomorphism} \label{BC_sec}

\subsection{Construction}

Again fix the depth-zero character $\chi$ on $T(F)_1$.  Let $(I,\rho)$ be the $\mathfrak s_\chi$-type described 
above.  Let $A^{F_r}$ denote the unique maximal $F_r$-split torus in $G$ which contains $A$, and note that $T = {\rm Cent}_G(A^{F_r})$.  Write $T_r$ for $T(F_r)$.  Consider $\chi_r := \chi \circ N_r$ as a depth-zero character on $T(F_r)_1$, and consider the corresponding inertial class ${\mathfrak s}_r := \mathfrak s_{\chi_r}$ for $G_r = G(F_r)$.  Let $(I_r, \rho_r)$ denote the $\mathfrak s_r$-type associated to the character $\chi_r$.  We denote the corresponding Hecke algebra (resp. its center) by $\mathcal H(G_r,\rho_r)$ (resp. $\mathcal Z(G_r,\rho_r)$).

There is a canonical morphism of algebraic varieties
\begin{align*}
N_r^* : \mathfrak X_{\mathfrak s} &\rightarrow \mathfrak X_{\mathfrak s_r} \\
(T,\xi)_G &\mapsto (T_r, \, \xi \circ N_r)_{G_r},
\end{align*}
where here $\xi$ denotes an extension of some $W(F)$-conjugate of $\chi$.  This induces an algebra homomorphism 
\begin{equation} 
N_r: \mathbb C[{\mathfrak X_{\mathfrak s_r}}] \rightarrow \mathbb C[{\mathfrak X_\mathfrak s}].
\end{equation}

\begin{defn} \label{bc_defn}
Define the {\em base change homomorphism} $b_r: \mathcal Z(G_r, \rho_r) \rightarrow \mathcal Z(G,\rho)$ to be the unique morphism making the following diagram commute:
$$
\xymatrix{
\mathbb C[{\mathfrak X_{\mathfrak  s_r}}] \ar[d]^{N_r} \ar[r]^{\beta}_{\sim} & \mathcal Z(G_r, \rho_r) \ar[d]^{b_r} \\
\mathbb C[{\mathfrak X_\mathfrak s}] \ar[r]^{\beta}_{\sim} & \mathcal Z(G,\rho).}
$$
\end{defn}

\subsection{Characterization by actions on principal series}

\begin{lemma}\label{left_ch_br_lem}  For any character $\xi: T(F)\rightarrow \mathbb C^\times$ which extends some $W(F)$-conjugate of $\chi$, define $\xi_r := \xi \circ N_r$, a character on $T(F_r)$ which extends some $W(F_r)$-conjugate of $\chi_r$.  Let $z \in \mathcal H(G_r, \rho_r)$.  Then $b_r(z)$ is the unique element in $\mathcal Z(G,\rho)$ which acts on every module $i^G_B(\xi)^\rho$ by the same scalar by which $z$ acts on $i^{G_r}_{B_r}(\xi_r)^{\rho_r}$.  In other words, 
\begin{equation} \label{left_ch_br}
\beta^{-1}(b_r(z))(\xi) = \beta^{-1}(z)(\xi_r).
\end{equation}
\end{lemma}
In terms of right actions this means: for $z \in \mathcal Z(G_r, \rho_r)$, $b_r(z)$ acts on the right on $i^G_B(\xi^{-1})^{\rho^{-1}}$ by the scalar by which $z$ acts on the right on $i^{G_r}_{B_r}(\xi_r^{-1})^{\rho_r^{-1}}$, that is,
\begin{equation} \label{ch-b_r}
ch_{\xi^{-1}}(b_r(z)) = ch_{\xi^{-1}_r}(z)
\end{equation}
(cf. Lemma \ref{ch_xi_lem}). 

This is reflective of the fact that $b_r$ is compatible with the involutions $\iota_r$ (resp. $\iota$) of $\mathcal H(G_r)$ (resp. $\mathcal H(G)$), in the sense that the following diagram commutes
$$
\xymatrix{
\mathcal Z(G_r,\rho_r) \ar[r]^{b_r} \ar[d]_{\iota_r} & \mathcal Z(G,\rho) \ar[d]_\iota \\
\mathcal Z(G_r ,\rho^{-1}_r) \ar[r]^{b_r} & \mathcal Z(G,\rho^{-1}).}
$$
Indeed, this follows from (\ref{beta/iota}) and the obvious compatibility of $N_r$ and $\iota : \mathbb C[{\mathfrak X_\mathfrak s}] \rightarrow \mathbb C[{\mathfrak X_\mathfrak s}]$.

\subsection{Compatibility with conjugation by $w \in W(F)$}

Fix $w \in W(F)$, and use the same symbol to denote its lift in $N_G(T) \cap K$.  The character $^w\chi$ is defined by $^w\chi(t) = \chi(w^{-1}tw)$. (Similarly define $^w\Phi(\cdot) = \Phi(w^{-1} \,\cdot\, w)$ for any suitable function $\Phi$.)  Also write $^wI := wIw^{-1}$ and $^wI^+ = wI^+w^{-1}$.   The character $^w\rho : \, ^wI/\, ^wI^+ \rightarrow \mathbb C^\times$ is defined using $^w\chi$; in fact $^w\rho(wiw^{-1}) = \rho(i)$ for $i \in I$.  There is an isomorphism $\mathcal H(G,I,\rho) ~ \widetilde{\rightarrow} ~ \mathcal H(G,\, ^wI,\, ^w\rho)$, given by $h \mapsto \, ^wh$. 

Let $\xi$ denote an extension of a $W(F)$-conjugate of $\chi$.  Write $^wB := wBw^{-1}$.  Then $\Psi \mapsto \, ^w\Psi$ gives an isomorphism
$$
i^G_{B}(\xi^{-1})^{\rho^{-1}} ~ \widetilde{\rightarrow} ~ i^G_{\,^wB}(\, ^w\xi^{-1})^{\, ^w\rho^{-1}}.
$$

This intertwines the right actions of $h \in \mathcal H(G,I,\rho)$ and $^wh \in \mathcal H(G,\ ^wI,\, ^w\rho)$, in the sense that
$$
^w(\Psi * h) = \, ^w\Psi * \, ^wh.
$$

Applying this to $h = z \in \mathcal Z(G,I,\rho)$ and taking Remark \ref{^wbeta} into account yields the next lemma.

\begin{lemma} \label{beta-ad(w)}
The following diagram is commutative:
$$
\xymatrix{
\mathbb C[{\mathfrak X_\chi}] \ar[r]^{\beta \,\,\,\,\,}_{\sim \,\,\,\,\,} \ar[d]^{=} & \mathcal Z(G,I,\rho) \ar[d]^{z \mapsto \, ^wz} \\
\mathbb C[{\mathfrak X_\chi}] \ar[r]^{\beta \,\,\,\,\,}_{\sim \,\,\,\,\,} & \mathcal Z(G,\, ^wI, \, ^w\rho).}
$$
\end{lemma}

The compatibility of $z \mapsto \, ^wz$ with base change follows.

\begin{lemma} \label{bc-ad(w)}
For any $w \in W(F)$, the following diagram is commutative:
$$
\xymatrix{
\mathcal Z(G_r, \rho_r) \ar[r]^{z \mapsto \, ^wz} \ar[d]^{b_r} & \mathcal Z(G_r, \, ^w\rho_r) \ar[d]^{b_r} \\
\mathcal Z(G, \rho) \ar[r]^{z \mapsto \, ^wz} & \mathcal Z(G, \, ^w\rho).}
$$
\end{lemma}

\section{Constant term homomorphism} \label{CT_sec}

\subsection{Abstract definition}

Let $M$ denote an $F$-Levi subgroup of $G$ which contains the torus $T$.  Let $\mathfrak s^M_\chi$ denote the inertial class associated to $\chi$ but for the group $M$ rather than the group $G$.  Let $\mathcal R_\chi(M)$ denote the corresponding principal series block of $\mathcal R(M)$, and let $\mathfrak X_{\chi}^M$ denote the corresponding variety of supercuspidal supports $(T,\xi)_M$.

There is a canonical surjective morphism of varieties
\begin{align} \label{ct*_def}
 c^{G*}_M : \mathfrak X^M_{\chi} &\rightarrow \mathfrak X_\chi \\
(T,\xi)_M &\mapsto (T,\xi)_G. \notag
\end{align}
This induces an algebra homomorphism
\begin{equation} \label{ct_def}
c^G_M: \mathbb C[{\mathfrak X_\chi}] \rightarrow \mathbb C[{\mathfrak X^M_\chi}].
\end{equation}
We call $c^G_M$ the {\em constant term homomorphism}.   In the next few subsections we will give a concrete description which will also justify this terminology.  We will also prove its compatibility with the base change homomorphism.  As the first step, note that the following diagram is obviously commutative:
\begin{equation} \label{ct-N_r}
\xymatrix{
\mathbb C[{\mathfrak X_{\chi_r}}] \ar[r]^{c^{G_r}_{M_r}} \ar[d]_{N_r} & \mathbb C[{\mathfrak X^{M_r}_{\chi_r}}] \ar[d]_{N_r} \\
\mathbb C[{\mathfrak X_\chi}] \ar[r]^{c^G_M} & \mathbb C[{\mathfrak X^M_\chi}].}
\end{equation}

\subsection{Refined Iwasawa decomposition}

Let us suppose $M$ is a {\em standard}\footnote{This assumption is convenient but not necessary for this subsection and the next.} $F$-Levi subgroup; this means that $M$ is a Levi factor of a {\em standard} (i.e., containing $B$) $F$-parabolic subgroup $P$  with unipotent radical $N$.  Fix a Levi decomposition $P = MN$.
 
Consider the Kostant representatives $W^P \subset W(F)$, i.e. the set of minimal coset representatives for the elements of $W_M(F)\backslash W(F)$, with respect to the Bruhat order on $W(F)$ defined by the Borel subgroup $B$.  We choose once and for all a lift $\dot{w} \in N_G(T)(F) \cap K$ for each $w \in W^P$; but from now on we denote $\dot{w}$ simply by $w$.  

Let $P_r$ denote $P(F_r)$, etc.  Let $I_r$ denote the Iwahori subgroup of $G_r := G(F_r)$ corresponding to $I$, so that $I_r \cap G(F) = I$.  Abbreviate the relative Weyl group $W(F_r)$ by $W_r$ and its set of Kostant representatives by $W_r^P$.  The refined Iwasawa decomposition states that 
\begin{equation} \label{ref_Iwasawa}
G_r = \coprod_{w \in W^P_r} P_r\,w\,I_r.
\end{equation}

\subsection{Levi factorization for $P \cap \, ^wI$}

\begin{lemma} \label{levi_fact}
For each $w \in W$, we have the factorization $P \cap \, ^wI = (M \cap \, ^wI)\,(N \cap \, ^wI).$
\end{lemma}

\begin{proof}
This follows from \cite{BT2}, 5.2.4 (cf. section \ref{corr_sec}).
\end{proof}

\subsection{The concrete definition of the constant term}

Fix an Iwahori subgroup $J$ which is $W(F)$-conjugate to $I$.

For a compactly-supported locally constant function $f$ on $G(F)$, we define the compactly-supported locally constant function $f^{(P)}$ on $M(F)$ by the formula
\begin{equation} \label{f^P_def}
f^{(P)}(m) = \delta_P^{1/2}(m) \int_{N(F)} f(mn) \, dn = \delta_P^{-1/2}(m) \int_{N(F)} f(nm) \, dn,
\end{equation}
where the Haar measure $dn$ on $N(F)$ is normalized to give $N(F) \cap J$ measure 1.
Warning: $f^{(P)}$ depends on the choice of $J$, and in what follows we will allow $J$ to vary.

Since $J/J^+ \cong T(F)_1/T(F)_1^+$, the character $\chi$ gives rise to a smooth character $\rho: J \rightarrow \mathbb C^\times$ as before, and $(J,\rho)$ is also a type for $\mathcal R_\chi(G)$.  Write the corresponding Hecke algebra (resp. its center) as $\mathcal H(G,J,\rho)$ (resp. $\mathcal Z(G,J,\rho)$) in order to emphasize the subgroup $J$. 

Let $J_M = J \cap M$.  Note that $J_M$ is an Iwahori subgroup of $M$ (cf. e.g. \cite{H09a}, Lemma 2.9.1).  We get the character $\rho: J_M \rightarrow \mathbb C^\times$ from $\chi$ and the isomorphism $J_M/J^+_M \cong T(F)_1/T(F)_1^+$, and $(J_M, \rho)$ is a type for $\mathcal R_\chi(M)$.

It is obvious that $f^{(P)} \in \mathcal H(M,J_M,\rho)$ if $f \in \mathcal H(G,J,\rho)$.  The next proposition is the main goal of this section.  

\begin{prop} \label{pres_center}  The operation $f \mapsto f^{(P)}$ sends $\mathcal Z(G,J,\rho)$ to $\mathcal Z(M,J_M,\rho)$, and the following diagram commutes:
$$
\xymatrix{
\mathbb C[{\mathfrak X_\chi}] \ar[d]^{c^G_M} \ar[r]^{\beta \,\,\,\,}_{\sim \,\,\,\,} & \mathcal Z(G,J,\rho) \ar[d]^{f \mapsto f^{(P)}} \\
\mathbb C[{\mathfrak X^M_\chi}] \ar[r]^{\beta \,\,\,\,}_{\sim \,\,\,\,} & \mathcal Z(M,J_M,\rho).}
$$
In particular, when restricted to $\mathcal Z(G,J,\rho)$ the operation $f \mapsto f^{(P)}$ is an algebra 
homomorphism and is independent of the choice of $P$ having $M$ as Levi factor. 
\end{prop} 
Because of this we will set $c^G_M(f) := f^{(P)}$ for $f \in \mathcal Z(G,J,\rho)$ and call it the {\em constant term of $f$}.

\subsection{Compatibility of base change and constant term homomorphisms}

Proposition \ref{pres_center} and equation (\ref{ct-N_r}) have the following corollary.

\begin{cor} \label{ct-bc_comp}
The following diagram commutes:
\begin{equation} \label{eq:ct-bc}
\xymatrix{
\mathcal Z(G_r,J_r, \rho_r) \ar[r]^{c^{G_r}_{M_r}} \ar[d]_{b_r} & \mathcal Z(M_r,J_{M_r},\rho_r) \ar[d]_{b_r} \\
\mathcal Z(G,J,\rho) \ar[r]^{c^G_M} & \mathcal Z(M,J_M,\rho).}
\end{equation}
\end{cor}

\subsection{Proof of Proposition \ref{pres_center}}

Fix an extension $\xi$ of $\chi$, and set $\xi_1 := \delta_P^{1/2} \xi$.  Consider the subspace 
$i^G_B(\xi^{-1}_1)^{\rho^{-1}} \subset i^G_B(\xi^{-1}_1)$ whose elements transform under the left $J$-action according to the character $\rho^{-1}$.   To prove Proposition \ref{pres_center}, we need to prove the following equivalent statement.

\begin{prop} \label{alt_pres_center}
Let $z \in \mathcal Z(G,J, \rho)$.  Then for every $\xi$, the element $z^{(P)}$ acts on the right on $i^M_{B_M}(\xi^{-1})^{\rho^{-1}}$ by $ch_{\xi^{-1}}(z)$, the scalar by which $z$ acts on the right on $i^G_B(\xi^{-1})^{\rho^{-1}}$. 
\end{prop}

\begin{lemma} \label{extn_lem} Given $\psi \in i^M_{B_M}(\xi^{-1})^{\rho^{-1}}$, there exists $\Psi \in i^G_B(\xi_1^{-1})^{\rho^{-1}}$ such that $\Psi|_M = \psi$.
\end{lemma}

\begin{proof}
We define $\Psi$ to be the unique element in $i^G_{B}(\xi^{-1}_1)^{\rho^{-1}}$ which is supported on $PJ$ and which extends $\psi$.  That is, for $mn \in MN$ and $j\in J$ we set 
$$\Psi(mnj) := \psi(m)\rho^{-1}(j).
$$
To see that $\Psi$ belongs to the desired space, we use the fact that $\delta_B^{1/2}\xi^{-1}_1 = \delta_{B_M}^{1/2}\xi^{-1}$.  To see it is well-defined, we use the Iwahori-factorization of $J$ with respect to $P = MN$, or more precisely a consequence of it: $J \cap P = J \cap M \cdot J \cap N$.   See Lemma \ref{levi_fact}.
\end{proof}

Proposition \ref{alt_pres_center} will be proved using the following lemma.  

\begin{lemma} \label{main_calc} Define the right convolution action $\ast$ of $\mathcal H(M,J_M, \rho)$ on $i^M_{B_M}(\xi^{-1})^{\rho^{-1}}$ using the measure $dm$ which gives $J_M$ measure 1.  Let $z \in \mathcal Z(G,J,\rho)$.  Let $\psi \in i^M_{B_M}(\xi^{-1})^{\rho^{-1}}$.  Then for any $y \in M$, we have the identity
\begin{equation} \label{eq:main_calc}
ch_{\xi^{-1}_1}(z) \cdot \, {\psi} (y) = \psi * [\delta_P^{1/2} z^{(P)}](y).
\end{equation}
\end{lemma}

\begin{proof}
In addition to $dm$ and $dn$ specified above, we let $dg$ (resp. $dj$) denote the Haar measure on $G$ (resp. $J$) giving $J$ measure 1.  Let $\Psi \in i^G_B(\xi_1^{-1})^{\rho^{-1}}$ denote the extension of $\psi$ from Lemma \ref{extn_lem}.  We have
\begin{align*}
\psi * [\delta_P^{1/2} \, z^{(P)}](y) &= \int_{M} \psi(m) \, \delta_P^{1/2}(m^{-1}y) \, z^{(P)}(m^{-1}y) \, dm \\
&= \int_M \int_N \Psi(mn) \, z((mn)^{-1}y) \, dn \, dm \\
&= \int_M \int_N \int_J \Psi(mnj) \, z((mnj)^{-1}y) \, dj \, dn \, dm.
\end{align*}
Now since $\Psi$ is supported on $MNJ$, we may use the substitution $g = mnj$ (cf. \cite{H09a}, $\S4.3$) to write this as
\begin{align*}
\int_G \Psi(g) \, z(g^{-1}y) \, dg &= (\Psi * z)(y) \\
&= ch_{\xi^{-1}_1}(z) \cdot \Psi(y).
\end{align*}
The lemma follows since $\Psi(y) = \psi(y)$.
\end{proof}

\smallskip

\noindent {\em Proof of Proposition \ref{alt_pres_center}}:  Recalling that $\xi^{-1} = \delta_P^{1/2} \xi^{-1}_1$, Lemma \ref{main_calc} shows that $\delta_P^{1/2} \, z^{(P)}$ acts on the right on $i^M_{B_M}(\delta_P^{1/2} \xi_1^{-1})^{\rho^{-1}}$ by $ch_{\xi^{-1}_1}(z)$, the scalar by which $z$ acts on the right on $i^G_B(\xi^{-1}_1)^{\rho^{-1}}$.  It follows (cf. \cite{H09a}, (4.7.3)) that $z^{(P)}$ acts on the right on $i^M_{B_M}(\xi_1^{-1})^{\rho^{-1}}$ by the same scalar $ch_{\xi^{-1}_1}(z)$.  Since $\xi_1$ ranges over all extensions of $\chi$ as $\xi$ does, this proves the proposition.  \qed

\section{Descent formulas}  \label{Des_sec}

\subsection{Preliminaries}

Fix once and for all Haar measures $dg$,$di$ on $G,I$ respectively, such that ${\rm vol}_{dg}(I) = {\rm vol}_{di}(I) = 1$.  On $P = MN$ we fix a (left) Haar measure $dp$ such that ${\rm vol}_{dp}(P \cap I) = 1$.  For $w \in W(F)$, let $dp_w$ (resp. $dm_w$, $dn_w$) denote the (left) Haar measure on $P$ (resp. $M$,$N$) such that ${\rm vol}_{dp_w}(P \cap \, ^wI)$ (resp. ${\rm vol}_{dm_w}(M \cap \, ^wI)$, ${\rm vol}_{dn_w}(N \cap \, ^wI)$) has the value 1.

Similar conventions will hold for measures on the groups $G_r, P_r, I_r$, etc.

Recall the Harish-Chandra function on $m \in M(F)$, defined by
$$
D_{G(F)/M(F)}(m) = {\rm det}\big( 1 - {\rm Ad}(m^{-1}); {\rm Lie}(G(F))/{\rm Lie}(M(F))\big).
$$

\subsection{Descent of (twisted) orbital integrals}

Let us recall the set-up: $\gamma \in G(F)$ is a semisimple element, $S$ denotes the $F$-split component of the center of $G_\gamma^\circ$, and $M = {\rm Cent}_G(S)$, an $F$-Levi subgroup.  We suppose $\gamma$ is not elliptic in $G(F)$, so that $M$ is a proper Levi subgroup of $G$.  Choose an $F$-parabolic $P = MN$ with $M$ as Levi factor.  We may assume $P$ and $M$ are standard.

We have $\gamma \in M(F)$ and $G^\circ_\gamma = M^\circ_\gamma$.  We assume $\gamma = \mathcal N \delta$ for an element $\delta \in G(F_r)$.  By Lemma 4.2.1 of \cite{H09a}, we may assume $\delta \in M(F_r)$.  

The twisted centralizer $G_{\delta \theta}$ of $\delta \theta$ is an inner form of $G_\gamma$ whose group of $F$-points is
$$
G_{\delta \theta}(F) = \{ g \in G(F_r) ~ | ~ g^{-1}\delta\theta(g) = \delta \}.
$$
We have $M^\circ_{\delta\theta} = G^\circ_{\delta \theta}$ (cf. \cite{H09a}, $\S4.4$).  We choose compatible measures on the inner forms $G^\circ_{\delta \theta}$ and $G^\circ_\gamma$, and use them to form the quotient measures $d\bar{g}$ in the (twisted) orbital integrals of $\phi \in C^\infty_c(G(F_r))$.  By definition,
\begin{equation} \label{TO_def}
{\rm TO}^{G_r}_{\delta\theta}(\phi) = \int_{G^\circ_{\delta \theta}\backslash G_r} \phi(g^{-1} \delta \theta(g)) \, d\bar{g}.
\end{equation}
Here and in what follows we denote $G^\circ_{\delta \theta}(F)$ simply by $G^\circ_{\delta \theta}$.

Note that $\rho_r\circ \theta = \rho_r$.  Using this and the refined Iwasawa decomposition (\ref{ref_Iwasawa}), the argument in \cite{H09a}, $\S4.4$ yields for $\phi \in \mathcal H(G_r, \rho_r)$ the descent formula
\begin{equation} \label{TO_descent}
{\rm TO}^{G_r}_{\delta \theta}(\phi) = |D_{G(F)/M(F)}(\gamma)|_F^{-1/2} \sum_{w \in W_r^P} {\rm TO}^{M_r}_{\delta \theta}\big((\,^{w,\theta}\phi)^{(P_r)}\big),
\end{equation}
where $^{w,\theta}\phi(g)   := \phi(w^{-1}g\theta(w))$.  Here it is understood that ${\rm TO}^{M_r}$ is formed using $dm_w$ and $(\cdot)^{(P_r)}$ is formed using $dn_w$.

Orbital integrals are special cases of twisted orbital integrals (take $r=1$ and $\delta = \gamma$).  For $\phi \in \mathcal Z(G_r, \rho_r)$, we have the following descent formula for $b_r(\phi) \in \mathcal Z(G,\rho)$:
\begin{equation} \label{O_descent1}
{\rm O}^G_\gamma(b_r(\phi)) = |D_{G(F)/M(F)}(\gamma)|_F^{-1/2} \sum_{w \in W^P} {\rm O}^M_\gamma[(\,^wb_r(\phi))^{(P)}].
\end{equation}
Using the compatibility of $b_r$ with conjugation by $w$ (Lemma \ref{bc-ad(w)}) and constant term (Corollary \ref{ct-bc_comp}, taking $J = \, ^wI$), we can write this as follows:
\begin{equation} \label{O_descent2}
{\rm O}^G_\gamma(b_r(\phi)) = |D_{G(F)/M(F)}(\gamma)|_F^{-1/2} \sum_{w \in W^P} {\rm O}^M_\gamma \big(b_r(( ^w\phi)^{(P_r)})\big).
\end{equation}

\subsection{Comparing descent formulas}

We compare the formulas (\ref{TO_descent}) and (\ref{O_descent2}) following the method of \cite{H09a}, $\S4.5$.  Note that $W^P = (W^P_r)^\theta$.  Thus the key point is the following lemma.  

\begin{lemma} \label{vanishing}
Assume $\phi \in \mathcal Z(G_r, \rho_r)$.  Then the summands in (\ref{TO_descent}) indexed by elements $w \in W^P_r$ with $\theta(w) \neq w$ are zero.
\end{lemma}

\begin{proof}
This is the analogue of Lemma 4.5.3 in \cite{H09a} which handled the case of parahoric Hecke algebras.  The proof in pp. 602-608 of loc.~cit. goes over to the present context nearly word-for-word; we omit the details.
\end{proof}

Comparing (\ref{TO_descent}) and (\ref{O_descent2}), it is now clear that the fundamental lemma for $\gamma, \delta, G$ follows from the fundamental lemma for $\gamma, \delta, M$.  Thus, by induction on the semisimple rank of $G$, we are reduced to considering elliptic semisimple elements $\gamma$.

\begin{Remark}
Assume $G_{\rm ad}$ is split over $F$.  Then $W^P = W_r^P$ and Lemma \ref{vanishing} is not needed to compare (\ref{TO_descent}) with (\ref{O_descent2}).
\end{Remark}

\section{Reductions} \label{reductions_section}

\subsection{Vanishing statement when $\gamma$ is not a norm}

\begin{lemma} \label{wlog_norm}
Let $\phi \in \mathcal Z(G_r, \rho_r)$.  If $\gamma$ is not a norm from $G_r$, then ${\rm SO}_\gamma(b_r(\phi)) = 0$.
\end{lemma}

We will prove the stronger fact that ${\rm O}_\gamma(b_r(\phi)) =0$.  The corresponding statement for spherical Hecke algebras was proved by Labesse \cite{Lab99}, Lemme 3.7.1.  The proof here is similar, but we 
give some details for the sake of completeness. 

First, our descent formula (\ref{O_descent2}) reduces us to the case where $\gamma$ is elliptic in $G$.  Consider the canonical map $p:G_{\rm sc} \rightarrow G$ and the abelian group
$$
{\bf H}^0_{ab}(F,G) = G(F)/p(G_{\rm sc}(F)).
$$
Proposition 2.5.3 of loc.~cit. shows that an elliptic element $\gamma$ is a norm from $G_r$ if and only if its image in ${\bf H}^0_{ab}(F,G)$ is a norm.  Now the required vanishing follows from the following result (cf. \cite{Lab99}, Lemme 3.7.1).

\begin{lemma}\label{3.7.1}
Let $f = b_r(\phi)$ for $\phi \in \mathcal Z(G_r,\rho_r)$.  Let $x \in G(F)$ be any element such that, for some character $\eta$ on the group ${\bf H}^0_{ab}(F,G)$ which is trivial on the norms, we have $\eta(x) \neq 1$.  Then $f(x) = 0$.
\end{lemma}

\begin{proof}
By pulling back along the projection $G(F) \rightarrow {\bf H}^0_{ab}(F,G)$, we view $\eta$ as a 
character $\eta: G(F) \rightarrow \mathbb C^\times$; thus the condition $\eta(x) \neq 1$ makes sense.

It is not hard to show that $\eta$ is necessarily an unramified character of $G(F)$, meaning it is trivial on $G(F)_1 = G(F) \cap G(L)_1$, where $G(L)_1$ is the kernel of the Kottwitz homomorphism $\kappa: G(L) \rightarrow X^*(Z(\widehat{G})^I)$ (cf. (\ref{K-hom})).  Equivalently, $\eta$ is trivial on every parahoric subgroup of $G(F)$ (as $G(F)_1$ is the group generated by the parahoric subgroups of $G(F)$, cf. \cite{HR1}).  The restriction of $\eta$ to $T(F)$ is also unramified (trivial on $T(F)_1$).  

Thus, $f\eta \in \mathcal H(G,\rho)$.   Now by examining the right convolution action of $f\eta$ on $i^G_B(\xi^{-1})^{\rho^{-1}}$, we see that $f\eta \in \mathcal Z(G,\rho)$ and in fact 
$$
ch_{\xi^{-1}}(f\eta) = ch_{(\eta \xi)^{-1}}(f).
$$
But by (\ref{ch-b_r}), this is 
$$
ch_{(\eta_r \xi_r)^{-1}}(\phi) = ch_{\xi_r^{-1}}(\phi) = ch_{\xi^{-1}}(f),
$$
the first equality holding since $\eta_r = {\rm triv}$.  But this implies that $f\eta = f$, and in particular
$$
f(x)(\eta(x)-1) = 0.
$$
Since $\eta(x)-1 \neq 0$, the lemma follows.
\end{proof}  
 
\subsection{Reductions when $\gamma$ is a norm: three lemmas}

Our strategy is as in \cite{H09a}, $\S5$ which handled the case of parahoric Hecke algebras.  The first task is to state three key lemmas, which are analogous to Lemmas 5.3.1, 5.3.2, and 5.3.3 from \cite{H09a}.  The proofs go along the same lines; we shall explain the objects where they differ from loc.~cit., but we shall omit detailed proofs.  We also use freely the notation from loc.~cit.

Choose a finite unramified extension $F' \supset F$ which contains $F_r$ and splits $G$.  Consider a $z$-extension of $F$-groups
$$
\xymatrix{
1 \ar[r] & Z \ar[r] & H \ar[r]^p & G \ar[r] & 1,}
$$
where $Z$ is a finite product of copies of $R_{F'/F}{\mathbb G}_m$.   As usual, we are assuming $H_{\rm der} = H_{\rm sc}$.  Recall that $p$ is surjective on $F_r$- and $F$-points.  Choose an extension of $\theta$ to an element, still denoted $\theta$, in ${\rm Gal}(F'/F)$.  As in loc.~cit., the norm homomorphism $N: Z(F_r) \rightarrow Z(F)$ is surjective and induces a surjective map $N: Z(F_r)_1 \rightarrow Z(F)_1$. 

Let $\lambda: Z(F) \rightarrow \mathbb C^\times$ denote a smooth character, and for $f \in C^\infty_c(H(F))$ set
$$
f_\lambda(h) = \int_{Z(F)} f(hz) \lambda^{-1}(z) \, dz
$$
where $dz$ is the Haar measure on $Z(F)$ giving $Z(F)_1$ measure 1.  Write $\lambda = 1$ for the trivial character, and $\overline{f}$ for the function $f_1$ when it is viewed as an element in $C^\infty_c(G(F))$.   Write $\lambda N$ for the character $\lambda \circ N: Z(F_r) \rightarrow \mathbb C^\times$. 

The depth-zero character $\chi$ on $T(F)_1$ determines a depth-zero character $\chi_H$ on $T_H(F)_1$, where $T_H := p^{-1}(T)$.  Denote by $I_H$ the Iwahori subgroup in $H$ corresponding to the Iwahori subgroup $I$ in $G$.  Let $\rho_H :I_H \rightarrow \mathbb C^\times$ denote the character constructed from $\chi_H$.

\begin{lemma} \label{5.3.1}
Let $\phi \in C^\infty_c(H(F_r))$ and $f \in C^\infty_c(H(F))$.  We have the following statements: 
\smallskip

\noindent {\rm (i)}  The functions $\phi, f$ are associated if and only if $\phi_{\lambda N},f_\lambda$ are associated for every $\lambda$.
\smallskip

\noindent {\rm (ii)}  Suppose that $\phi \in \mathcal H(G_r,\widetilde{K}, \widetilde{\rho})$ (resp. $f \in \mathcal H(G,K,\rho)$) for a compact open subgroup $\widetilde{K} \subset G_r$ and character $\widetilde{\rho}: \widetilde{K} 
\rightarrow \mathbb C^\times$ (resp. $K \subset G$ and $\rho: K \rightarrow \mathbb C^\times$) such that
\begin{itemize}
\item $N(\widetilde{K} \cap Z(F_r)) = K \cap Z(F)$,
\item $\widetilde{K} \supset (1-\theta)(Z(F_r))$, and
\item $\widetilde{\rho}|_{\widetilde{K} \cap Z(F_r)} = \rho \circ N|_{\widetilde{K} \cap Z(F_r)}$.
\end{itemize}
Then in (i) we only need to consider characters $\lambda$ such that 
$$
\lambda|_{K \cap Z(F)} = \rho^{-1}|_{K \cap Z(F)}.
$$
\smallskip

\noindent {\rm (iii)}  If $\phi \in \mathcal H(H_r, \rho_{Hr})$ and $f \in \mathcal H(H,\rho_H)$, then in (i) we need consider only the characters $\lambda$ with $\lambda|_{Z(F)_1} = \chi^{-1}|_{Z(F)_1}$.  
\smallskip

\noindent {\rm (iv)}  The pair $\phi_1,f_1$ are associated if and only if $\overline{\phi}, \overline{f}$ are associated.
\end{lemma}  

\begin{proof}
Parts (i) and (iv) are proved in \cite{H09a}.  Part (ii) is proved that same way as loc.~cit. Lemma 5.3.1(ii).  Part (iii) follows from part (ii), taking $K = I_H$, $\widetilde{K} = I_{Hr}$ and $\widetilde{\rho}= \rho_{Hr}$.  
\end{proof}

The following lemma is proved similarly to \cite{H09a}, Lemma 5.3.2.

\begin{lemma} \label{5.3.2}
The map $\phi \mapsto \overline{\phi}$ determines a surjective homomorphism $\mathcal Z(H_r, \rho_{Hr}) \rightarrow \mathcal Z(G_r,\rho_r)$.  It is compatible with the base change homomorphisms in the sense that
\begin{equation} \label{b_v_bar}
b_r(\overline{\phi}) = \overline{b_r(\phi)}.
\end{equation}
\end{lemma}

Now for any character $\lambda: Z(F)\rightarrow \mathbb C^\times$ such that $\lambda|_{Z(F)_1} = \chi^{-1}|_{Z(F)_1}$, consider the algebra $\mathcal H_{\rho_{Hr}, \lambda N}$ which consists of the locally constant functions $\phi$ on $H(F_r)$ which are compactly supported modulo $Z(F_r)$ and which transform under right and left multiplications by $i \in I_{Hr}$ (resp. $z \in Z(F_r)$) by $\rho^{-1}_{Hr}(i)$ (resp. $\lambda N(z)$).  One proves the following lemma by imitating the argument of \cite{H09a}, Lemma 5.3.3, using this algebra in place of the algebra $\mathcal H_{J,\chi N}$ used there.  Along the way, the argument relies on our Lemma \ref{5.3.1}(iii-iv) and Lemma \ref{5.3.2}.

\begin{lemma} \label{5.3.3}
Assume $Z = Z(H)$, so that $G = H_{\rm ad}$.  Suppose that $\overline{\phi}, b_r(\overline{\phi})$ are associated for every $\phi \in \mathcal Z(H_r, \rho_{Hr})$.  Then $\phi, b_r(\phi)$ are associated for every $\phi \in \mathcal Z(H_r, \rho_{Hr})$.  
\end{lemma}

\subsection{The reduction steps}

We may assume $\gamma = \mathcal N(\delta)$.  We follow \cite{H09a}, $\S5.4$.

\smallskip

(1)  {\em We may assume $G_{\rm der} = G_{\rm sc}$.}  Using Lemmas \ref{5.3.1} and \ref{5.3.2}, we see that the fundamental lemma for $H$ implies the fundamental lemma for $G$.

\smallskip

(2) {\em We may assume $\gamma$ is elliptic.}  We use our descent formulas (\ref{TO_descent}) and (\ref{O_descent2}), along with Lemma 4.2.1 of \cite{H09a}.  (See \cite{H09a}, $\S5.4(2)$.)

\smallskip

(3) {\em We may assume $\gamma$ is regular.}  This is proved by Clozel \cite{Cl90}, Prop. 7.2.

\smallskip

(4) {\em We may assume $G$ is such that $G_{\rm der} = G_{\rm sc}$ and $Z(G)$ is an induced torus (of the same form as $Z$ above).}  We work under assumptions (1-3), and follow Clozel's argument in \cite{Cl90}, section 6.1(b), with a few modifications.  Where Clozel uses the Satake isomorphism, we use the inverse $\beta^{-1}$ of the Bernstein isomorphism (\ref{isom1}).  The homomorphism $j$ used by Clozel is replaced by the natural homomorphism
$$
j: \mathcal Z(G,\rho) \hookrightarrow \mathcal Z(G', \rho')
$$
coming from the exact sequence 
$$
1 \rightarrow  G \rightarrow  G' \rightarrow  Q \rightarrow  1
$$
of loc.~cit. in the obvious way.  Here $G$ and $G'$ are unramified groups with the properties that $G_{\rm der} = G_{\rm sc}$, $G'_{\rm der} = G'_{\rm sc}$, and $Z(G')$ is an induced torus; we assume the fundamental lemma holds for the group $G'$ and our present goal is to deduce that it must hold for $G$.  The main point is to justify the analogues of Clozel's equalities (which we write here in slightly different notation) 
\begin{align} \label{Clozel_c=c'}
{\rm SO}^G_\gamma(f,dt,dg) &= c_T \, {\rm SO}^{G'}_\gamma(jf, dt', dg') \\
{\rm SO}^{G_r}_{\delta \theta}(\phi,dt,dg_{F_r}) &= C_T \, {\rm SO}^{G'_r}_{\delta \theta}(j\phi, dt', dg'_{F_r}).  \notag
\end{align}
We need to prove the equalities (\ref{Clozel_c=c'}) for all functions $f \in \mathcal Z(G,\rho)$ and $\phi \in \mathcal Z(G_r,\rho_r)$ and for all ($\theta$-)regular ($\theta$-)elliptic elements $\gamma$ and $\delta$ in $G$ whose ($\theta$-)centralizer is the elliptic $F$-torus $T$\footnote{Warning: in this subsection and the next, $T$ will denote an arbitrary elliptic torus, not the standard Cartan we denote with this symbol elsewhere!}.  Crucially, the constants $c_T$ and $C_{T}$ making (\ref{Clozel_c=c'}) hold must depend only on the torus $T$ and various choices of measures ($dt$, $dt'$, etc.), and we must have $c_T = C_T$.

The proof of the existence of the constant $C_T$ is unfortunately rather laborious.  The details we present apply as well to the spherical case considered by Clozel and to other settings, so we hope their potential use to the reader justifies their inclusion here.

To prove the existence of the constants $c_T$ and $C_T$, one first examines those functions $f,\phi$ which correspond under Morris's Hecke algebra isomorphisms \cite{Mor} to Bernstein functions in suitable Iwahori-Hecke algebras, and uses the fact that the homomorphisms $j$ are compatible with the Hecke algebra isomorphisms.  Let us make this more precise, using some notation from \cite{HR2}, $\S9$, and \cite{Ro}.  The Hecke algebra isomorphisms take the form 
\begin{equation} \label{HA_isom}
\mathcal H(G,\rho) \cong H_{\rm aff}(\widetilde{W}_H) ~\widetilde{\otimes} ~ \mathbb C[W_\chi/W^\circ_\chi],
\end{equation}
where $H$ is a certain connected reductive $F$-group with finite Weyl group $W^\circ_\chi$, with extended affine Weyl group $\widetilde{W}_H = X_*(A) \rtimes W^\circ_\chi$, and with associated affine Hecke algebra $H_{\rm aff}(\widetilde{W}_H)$.  Each $W^\circ_\chi$-orbit representative $\lambda \in X_*(A)$ gives rise to a Bernstein function $z_{\lambda} \in Z(H_{\rm aff}(\widetilde{W}_H))$ and for each $W_\chi$-orbit representative $\mu$ we have the Bernstein function
$$
Z_\mu := \sum_{\lambda} z_{\lambda}
$$
where $\lambda$ ranges over $W^\circ_\chi$-orbit representatives in $W_\chi \cdot \mu$;  the $Z_\mu$ form a basis for the center of the right hand side of (\ref{HA_isom}).  If $f \in \mathcal Z(G,\rho)$ corresponds to $Z_\mu$, then $jf \in \mathcal Z(G',\rho')$ corresponds to $Z'_\mu$ defined in the same way as $Z_\mu$, and we may write
\begin{align*}
f &= \sum_{w \in \mathcal S} a_w ~ [In_w I]_{\breve \chi} \\
jf &= \sum_{w \in \mathcal S'} a'_w ~ [I' n_w I']_{\breve \chi'},
\end{align*}
for constants $a_w, a'_w \in \mathbb C$ indexed by certain finite subsets $\mathcal S, S' \subset \widetilde{W}_H$ determined by $\mu$.  (We are writing the basis elements $[I n_w I]_{\breve \chi} \in \mathcal H(G,\rho)$, etc., following the notation of \cite{HR2},$\S9$.)  In fact $\mathcal S = \mathcal S'$ and the corresponding coefficients are equal, namely $a_w = a'_w$.  This follows, in view of the precise form of the Hecke algebra isomorphisms (see \cite{HR2}, Theorem 9.3.1 for the case where $G$ is split and $W^\circ_\chi = W_\chi$), from the 
corresponding equality of Iwahori-Matsumoto coefficients of $Z_\mu$ and $Z'_\mu$ (which is easy to prove from 
their definitions).   

In light of this, the existence of the constants $c_T$, $C_T$ follows via the next lemma.  We state this only in the twisted case (the untwisted case being the special case $r=1$ of the twisted case).  Recall that $\widetilde{W}_r$ denotes the extended affine Weyl group for $G_r$.

\begin{lemma} \label{C_T_lem}
There exists a constant $C_T$ such that, for every $\delta \in G_r$ with elliptic regular norm in $T$, and for every $w \in \widetilde{W}_r$, we have 
\begin{equation} \label{C_T:eq}
{\rm SO}^{G_r}_{\delta \theta}([I_r n_w I_r]_{\breve \chi_r}) = C_T \, {\rm SO}^{G'_r}_{\delta \theta}([I'_r n_w I'_r]_{\breve \chi'_r}).
\end{equation}
\end{lemma}

We postpone the proof of this lemma to the following subsection.

The matching of $jf, j\phi$ will force that of $f,\phi$ provided that $c_T =C_{T}$ for all maximal tori $T$ in $G$.  Since the homomorphisms $j$ are compatible with the base change homomorphisms, this will show that the fundamental lemma for $G'$ indeed implies the fundamental lemma for $G$.  In his situation, Clozel justifies the equality $c_T= C_T$ by invoking the matching of unit elements $1_K$ in spherical Hecke algebras and using the fact that orbital integrals of $1_K$ are never identically zero on any torus $T$ in $G$.

In our case, the construction of $c_T, C_T$ given in the proof of Lemma \ref{C_T_lem} shows that they do not depend on $\chi$ and they may be understood by taking $\chi = {\rm triv}$ and $w = 1$.  Thus, our equality $c_T = C_T$ follows from the fact that $1_{I_r}$ and $1_I$ are associated \cite{Ko86b}, and from the fact that the stable orbital integrals of $1_I$ do not vanish identically on any torus $T$ in $G$.  This last fact is evident for elliptic tori by the description of ${\rm O}_\gamma(1_I)$ in terms of the cardinality of the (non-empty) fixed-point set for the action of an elliptic $\gamma \in G$ on the building for $G$.  (Although we do not need this, the case of non-elliptic tori follows, by using our descent formula (\ref{O_descent1}), with $b(\phi_r)$ replaced by $1_I$, to reduce to the elliptic case.)

\smallskip

(5) {\em We may assume $G$ is adjoint.}  We know by (4) that the fundamental lemma holds provided it holds for all groups having the form of $H$.  By Lemma \ref{5.3.3}, the fundamental lemma for $G := H_{\rm ad}$ implies the fundamental lemma for $H$.  

\smallskip

(6)  {\em We may assume $\gamma$ is strongly regular elliptic.}  This is explained by Clozel \cite{Cl90}, p. 292 (see \cite{H09a}, $\S5.4(6)$).

\medskip

\noindent {\em Conclusion:}  We may assume that $G = G_{\rm ad}$ and $\gamma$ is a 
strongly regular elliptic semisimple element.  


\subsection{Proof of Lemma \ref{C_T_lem}}

I am very grateful to Robert Kottwitz for his invaluable help with this lemma.   Recall that the $F$-unramified groups $G,G'$ we are concerned with satisfy $G_{\rm sc} = G_{\rm der} = G'_{\rm der} = G'_{\rm sc}$ and $G/Z = G'/Z'$ where $Z$ resp. $Z'$ denotes the center of $G$ resp. $G'$\footnote{Much of the following discussion does not require $G_{\rm der}$ to be simply-connected.}.  By construction, the group $G'$ sits in a commutative diagram
$$
\xymatrix{
1 \ar[r] & G \ar[r] & G' \ar[r]^{p} & Q \ar[r] & 1 \\
1 \ar[r] & S \ar@{^(->}[u] \ar[r] & S' \ar@{^(->}[u] \ar[r] & Q \ar@{=}[u] \ar[r] & 1}
$$
and the (maximal Cartan) tori $S, \, S'$ and $Q$ are all unramified over $F$.  

Let $A_G$ denote the maximal $F$-split torus in the center of $G$.  Set $\Lambda = X_*(S)/Q^\vee$ where $Q^\vee$ denotes the coroot lattice.  Write $\Lambda(F)$ for the ${\rm Gal}(\bar{F}/F)$-fixed points on $\Lambda$.  Recall from \cite{Ko97} that there is a {\em surjective} homomorphism
\begin{equation} \label{kappa_surj}
\kappa: G(F) \twoheadrightarrow \Lambda(F)
\end{equation}
(cf. (\ref{K-hom})).  Let $G(F)_1$ denote the kernel of $\kappa$.  By Lemma \ref{surj_lem}(iii) below, the canonical map $S(F)_1 \rightarrow (G/G_{\rm der})(F)_1$ is surjective and hence $G(F)_1 = G_{\rm sc}(F) \, S(F)_1$.  Since $S(F)_1 \subset I$ (i.e. we may assume $S$ and $I$ are in good relative position), we have
\begin{equation} \label{G_scI}
G(F)_1 = G_{\rm sc}(F) \, I.
\end{equation}
Let $\Lambda_{\rm ad}$ denote the analogue of $\Lambda$ for the group $G_{\rm ad}$.
There is a short exact sequence
$$
0 \rightarrow X_*(A_G) \rightarrow \Lambda(F) \rightarrow \Lambda_{\rm ad}(F),
$$
and hence $\Lambda(F)/X_*(A_G)$ is a finite abelian group.

We can apply the constructions above to get corresponding objects $A_{G'}$, $\Lambda'$, etc., for $G'$.  There is a canonical map $\Lambda(F)/X_*(A_G) \rightarrow \Lambda'(F)/X_*(A_{G'})$.  It turns out that this is injective (one ingredient is the fact that $Q = Z'/Z$, which comes from the construction of $G'$ as in \cite{Cl90}).  Also, let $I'_r$ denote the Iwahori subgroup of $G'(F_r)$ corresponding to 
$I_r \subset G(F_r)$.  

We will need the following lemmas.  For the rest of this subsection, we assume all elements $\delta_0, \delta$, etc. belonging to $G(F_r)$ or $G'(F_r)$ have strongly regular semisimple norms.  Let $D := G/G_{\rm der}$ resp. $D' = G'/G'_{\rm der}$ and denote by $c: G' \rightarrow D'$ the projection homomorphism.

\begin{lemma} \label{surj_lem}  The following statements hold:

\noindent \item[(i)]  $p: S'(F)_1 \rightarrow Q(F)_1$ is surjective.

\noindent \item[(ii)] $p(G'(F)_1) = Q(F)_1$.

\noindent \item[(iii)] If $H \subset G$ is any maximal $F$-torus, then $c: H(F)_1 \rightarrow D(F)_1$ is surjective.

\noindent \item[(iv)] If $\delta \in G'(F_r)$, we have $G'(F)_1 \subseteq G'_{\delta \theta}(F)_1 \, G_{\rm sc}(F_r)$.  
\end{lemma}

\begin{proof}
Part (i).  Since $S$, $S'$, and $Q$ are all unramified tori, their Kottwitz homomorphisms yield a commutative diagram with exact rows and surjective vertical arrows
$$
\xymatrix{
1 \ar[r] & S(L) \ar[r] \ar[d] & S'(L) \ar[r]^p \ar[d] & Q(L) \ar[r] \ar[d] & 1 \\
0 \ar[r] & X_*(S) \ar[r] & X_*(S') \ar[r]^p & X_*(Q) \ar[r] & 0.}
$$
The snake lemma gives an exact sequence
$$
\xymatrix{
1 \ar[r] & S(L)_1 \ar[r] & S'(L)_1 \ar[r]^p & Q(L)_1 \ar[r] & 1.}
$$
Now (i) follows since $H^1(\langle \sigma \rangle, S(L)_1)=0$ (\cite{Ko97},(7.6.1)).  For part (iii), set $H_{\rm sc} := H \cap G_{\rm sc}$, and argue the same way; we get the exact sequence
$$
\xymatrix{
1 \ar[r] & H_{\rm sc}(L)_1 \ar[r] & H(L)_1 \ar[r]^c & D(L)_1 \ar[r] & 1}
$$
because $X_*(H_{\rm sc})_I$ is torsion-free.  

Now, part (ii) follows from part (i), since (i) gives the non-trivial inclusion in the chain
$$
Q(F)_1 \subseteq p(S'(F)_1) \subseteq p(G'(F)_1) \subseteq Q(F)_1.
$$

Finally, we deduce (iv) from (iii) applied to the group $G'$ with $H := G'_\gamma$, where $\gamma \in G'(F)$ is the norm of $\delta \in G'(F_r)$.  The groups $G'_{\delta \theta}$ and $G'_\gamma$ are inner forms over $F$ hence are $F$-isomorphic since they are tori.  There is an explicit isomorphism over $F_r$ (\cite{Ko82}), yielding a commutative diagram
$$
\xymatrix{
G'_{\delta \theta}(F_r)_1 \ar[d]^{\wr} \ar[r]^c & D'_{c(\delta) \theta}(F_r)_1 \ar[d]^{\wr} \\
 G'_\gamma(F_r)_1 \ar[r]^c & D'(F_r)_1.}
$$
Because of (iii), we see that the horizontal arrows are surjective and remain so upon taking ${\rm Gal}(F_r/F)$-invariants.   Doing so to the top arrow just gives the natural map
$$
c: G'_{\delta \theta}(F)_1 \rightarrow D'(F)_1,
$$
hence this is surjective.  Since ${\rm ker}(c) = G'_{\rm sc} = G_{\rm sc}$, (iv) follows easily.
\end{proof}

We use the symbol $w$ to denote any lift $n_w \in N_G(S)(F_r)$ for $w \in \widetilde{W}_r$, which we identify with $N_G(S)(F_r)/S(F_r)_1$.

\begin{lemma} \label{relevant_lem}
Let $\delta_0 \in G(F_r)$.  Suppose $\delta' \in G'(F_r)$ is stably $\theta$-conjugate to $\delta_0$ and belongs to $I'_r w I'_r$.  Then $\delta'$ is $I'_r$-$\theta$-conjugate to an element in $G(F_r) \cap I'_r w I'_r$.
\end{lemma}

\begin{proof}
First note that since $p(I'_r) \subseteq Q(F_r)_1$ and $p(w) = 1$, we have $p(\delta') \in Q(F_r)_1$.   
Writing $g'^{-1} \delta_0 \theta(g') = \delta'$ for some $g' \in G'(\bar{F})$ and applying $p$, one checks 
that $N_r(p(\delta'))=1$.  Now since $H^1(F_r/F, Q(F_r)_1) = 1$ (comp. (\ref{second_H1=0})), it follows that $p(\delta') \in (1-\theta)(Q(F_r)_1)$.  By Lemma \ref{surj_lem} (i), $p: S'(F_r)_1 \rightarrow Q(F_r)_1$ is surjective.  Hence there exists $s' \in S'(F_r)_1$ such that $p(s'\theta(s')^{-1}) = p(\delta')$.  But then $s'^{-1} \delta' \theta(s') \in G(F_r)$. 
\end{proof}

Now we fix an element $\delta_0 \in G(F_r)$ which we assume is $\theta$-elliptic and $\theta$-regular 
semisimple.  Also fix $w \in \widetilde{W}_r$ as above, and abbreviate the functions $[I_r n_w I_r]_{\breve \chi_r}$ resp. $[I'_r n_w I'_r]_{\breve \chi'_r}$ by $1_w$ resp. $1'_w$.  We seek to compare ${\rm SO}^{G_r}_{\delta_0 \theta}(1_w)$ with ${\rm SO}^{G'_r}_{\delta_0 \theta}(1'_w)$.  

Lemma \ref{relevant_lem} implies that if a $G'(F_r)$-$\theta$-conjugacy class in the stable $G'$-$\theta$-conjugacy class of $\delta_0$ contributes to ${\rm SO}_{\delta_0\theta}^{G'_r}(1'_w)$, then that class meets $G(F_r)$.  Hence we may write ${\rm SO}^{G'_r}_{\delta_0\theta}(1'_w)$ as a sum of terms ${\rm TO}^{G'_r}_{\delta \theta}(1'_w)$ where $\delta \in G(F_r)$ is stably $G'$-$\theta$-conjugate to $\delta_0$.

Now we fix such a $\delta$.  We will express ${\rm TO}^{G'_r}_{\delta \theta}(1'_w)$ in terms of certain integrals ${\rm TO}^{G_r}_{\delta_\lambda \theta}(1_w)$.  However, we start by studying ${\rm TO}^{G_r}_{\delta \theta}(1_w)$. Since $\delta$ is $\theta$-elliptic, we may use the following version of twisted orbital integrals
$$
{\rm TO}^{G_r}_{\delta \theta}(1_w) = \int_{A_G(F)\backslash G(F_r)} 1_w(g^{-1}\delta \theta(g)) \, \frac{dg}{dz},
$$
where ${\rm vol}_{dg}(I_r) = 1$ and ${\rm vol}_{dz}(A_G(\mathcal O)) = 1$.  This integral vanishes unless some $\theta$-conjugate $\delta_1$ of $\delta$ satisfies $\kappa(\delta_1) = \kappa(w)$.  So from now on we may assume that $\kappa(\delta) = \kappa(w)$.   We embed $X_*(A_G)$ in $A_G(F) \subset G(F)$ by $\nu \mapsto \varpi^\nu$ for any uniformizer $\varpi$.  The integral can be written as 
$$
\sum_{g \in X_*(A_G)\backslash G(F_r)/I_r} 1_w(g^{-1} \delta\theta(g)).
$$
Clearly the sum may be taken only over elements such that $g^{-1}\delta\theta(g) \in I_r n_w I_r$.  For these, $\kappa(g) \in \Lambda(F_r)/X_*(A_G)$ is fixed by $\theta$, i.e., $\kappa(g) \in \Lambda(F)/X_*(A_G)$.   We can then write the above sum as
$$
\sum_{\lambda \in \Lambda(F)/X_*(A_G)} ~ \sum_{\overset{g \in X_*(A_G) \backslash G_r/I_r}{\kappa(g) = \lambda}} 1_w(g^{-1}\delta \theta(g)).
$$
For each $\lambda \in \Lambda(F)$, use (\ref{kappa_surj}) to {\em choose} a $g_\lambda \in G(F)$ with $\kappa(g_\lambda) = \lambda$.
Using (\ref{G_scI}), we see that
\begin{equation*}
\{ g \in G(F_r)/I_r ~ | ~  \kappa(g) = \lambda \} = g_\lambda G_{\rm sc}(F_r)I_r/I_{r} 
= g_\lambda G_{\rm sc}(F_r)/I_{{\rm sc},r},
\end{equation*}
where $I_{{\rm sc},r} := I_r \cap G_{\rm sc}(F_r)$.  Set $\delta_\lambda := g^{-1}_\lambda \delta g_\lambda$.  We get the expression
\begin{equation} \label{TO^G_lambda}
{\rm TO}^{G_r}_{\delta \theta}(1_w) = \sum_{\lambda \in \Lambda(F)/X_*(A_G)} ~ \sum_{g \in G_{\rm sc}(F_r)/I_{{\rm sc},r}} 1_w(g^{-1} \delta_\lambda \theta(g)).
\end{equation}

Now we apply the same reasoning to $G'$.  Since $G'_{\rm sc} = G_{\rm sc}$, we find
\begin{equation} \label{TO^G'_lambda}
{\rm TO}^{G'_r}_{\delta \theta}(1'_w) = \sum_{\lambda' \in \Lambda'(F)/X_*(A_{G'})} ~ \sum_{g \in G_{\rm sc}(F_r)/I_{{\rm sc},r}} 1_w(g^{-1} \delta_{\lambda'} \theta(g)).
\end{equation}
(We used $1'_w(g^{-1}\delta_{\lambda'}\theta(g))= 1_w(g^{-1}\delta_{\lambda'}\theta(g))$ here.)   Now it follows that
\begin{equation} \label{TO_final}
{\rm TO}^{G'_r}_{\delta\theta}(1'_w) = \sum_{\lambda' \in \frac{\Lambda'(F)/X_*(A_{G'})}{\Lambda(F)/X_*(A_G)}} {\rm TO}^{G_r}_{\delta_{\lambda'}\theta}(1_w).
\end{equation}

By Lemma \ref{relevant_lem}, all relevant $G'(F_r)$-$\theta$-conjugacy classes in the stable $\theta$-conjugacy class of $\delta_0$ can be represented by elements $\delta$ in $G(F_r)$, so in forming ${\rm SO}^{G'_r}_{\delta_0\theta}(1'_w)$ we only sum ${\rm TO}^{G'_r}_{\delta \theta}(1'_w)$ over such elements $\delta \in G(F_r)$, for which the relation (\ref{TO_final}) holds.  Each such $\delta$ gives rise to elements $\delta_{\lambda'} = g_{\lambda'}^{-1} \delta g_{\lambda'}$ for various $g_{\lambda'} \in G'(F)$.  It follows that $\delta_{\lambda'} \in G(F_r)$ and is stably $\theta$-conjugate to $\delta$, hence to $\delta_0$.  

The following helps us understand the right hand side of (\ref{TO_final}).
\medskip

\begin{lemma} \label{saved} We have the following statements.

\noindent {\rm (1)}  Two elements of the form $\delta_{\lambda'}$ and $\delta_{\lambda''}$ are in the same $G(F_r)$-$\theta$-conjugacy class if and only if $\lambda'$ and $\lambda''$ belong to the same coset in $\frac{\Lambda'(F)/X_*(A_{G'})}{\Lambda(F)/X_*(A_G)}$.

\noindent {\rm (2)}  Any element $\delta_1 \in G(F_r) \cap I'_r w I'_r$ belonging to the $G'(F_r)$-$\theta$-conjugacy class of $\delta$ is $G(F_r)$-$\theta$-conjugate to some element of the form $\delta_{\lambda'}$.

\end{lemma}

\begin{proof} It is clearly enough to prove (1) in the special case $\lambda'' = 0$.  First assume there exists $g \in G(F_r)$ such that $\delta = (g_{\lambda'} g)^{-1}\delta \theta(g_{\lambda'} g)$.  Then $g_{\lambda'} g \in G'_{\delta \theta}(F)$.  Write $T'$ for the $F$-torus $G'_{\delta\theta}$ and $\kappa_{T'} : T'(F) \rightarrow X_*(T')^{{\rm Gal}(\bar{F}/F)}$ for the corresponding Kottwitz homomorphism.  We see that $\kappa_{T'}(g_{\lambda'}g) \in X_*(T')^{{\rm Gal}(\bar{F}/F)} = X_*(A_{G'})$, the equality holding since $T'$ is elliptic over $F$.  Recall that $\kappa(g_{\lambda'}g)$ is the image of $\kappa_{T'}(g_{\lambda'}g)$ in  $\Lambda'(F)$.  Thus $\lambda' \in X_*(A_{G'}) \Lambda(F_r) \cap \Lambda'(F) = X_*(A_{G'}) \Lambda(F)$.  

For the converse, suppose that $\lambda' \in X_*(A_{G'}) \, \Lambda(F)$.  Note that $X_*(A_{G'})$ can be embedded into $A_{G'}(F) \subset G'_{\delta \theta}(F)$ by $\nu \mapsto \varpi^\nu$, and that $\kappa_{G'}(\varpi^\nu)$ is the image of $\nu$ in $\Lambda'(F)$.  Then using Lemma \ref{surj_lem}(iv), we have
\begin{align*}
g_{\lambda'} &\in A_{G'}(F) \, G(F) \, G'(F)_1 \\
&\subseteq G'_{\delta \theta}(F) \, G(F_r).
\end{align*}
This implies that $\delta_{\lambda'} = g^{-1}_{\lambda'} \delta \theta(g_{\lambda'})$ is $G(F_r)$-$\theta$-conjugate to $\delta$.  

For (2), assume $\delta_1 = g'^{-1}\delta \theta(g')$ for some $g' \in G'(F_r)$.  Since $\kappa(\delta_1) = \kappa(w) = \kappa(\delta) \in \Lambda(F_r)$, we see $\kappa(g') \in \Lambda'(F)$; suppose $\kappa(g') = \lambda'$.  Then $g' \in g_{\lambda'} G'(F_r)_1 \subseteq g_{\lambda'} G_{\rm sc}(F_r) \, I'_r$ (by (\ref{G_scI})).  We can write
\begin{equation} \label{above}
i'^{-1}\delta_1 \theta(i') = h^{-1}\delta_{\lambda'} \theta(h)
\end{equation}
for some $i' \in I'_r$ and $h \in G_{\rm sc}(F_r)$.  

Now $\delta_1 \in G(F_r)$ means that $p(\delta_1) = 1$.  Also, the right hand side of (\ref{above}) is trivial under $p$.  It follows that $p(i') \in Q(F_r)_1^\theta = Q(F)_1$.  Now we apply Lemma \ref{surj_lem}(ii,iv) to $\delta_{1}$.  We see that $p(i') \in Q(F)_1 \subseteq p(G'(F)_1) \subseteq p(G'_{\delta_1 \theta}(F)_1)$ (since $p(G_{\rm sc}(F_r))=1$).  This means that
$$
i' \in G'_{\delta_1\theta}(F)_1 \, G(F_r).
$$
Thus the left hand side of (\ref{above}) is $G(F_r)$-$\theta$-conjugate to $\delta_1$, and it follows that $\delta_1$ is $G(F_r)$-$\theta$-conjugate to $\delta_{\lambda'}$.
\end{proof}
 
Finally, Lemma \ref{relevant_lem}, (\ref{TO_final}), and Lemma \ref{saved} together imply that
$$
{\rm SO}^{G'_r}_{\delta_0\theta}(1'_w) = {\rm SO}^{G_r}_{\delta_0\theta}(1_w)
$$
with the normalizations for twisted orbital integrals we fixed above.  This completes the proof of Lemma \ref{C_T_lem}.
\qed

\section{Labesse's elementary functions adapted to $\mathcal R_{\chi_r}(G_r)$} \label{elem_fcns_sec}

\subsection{Definition}

We will now construct the analogues of Labesse's (twisted) elementary functions \cite{Lab90} which are adapted to the Bernstein component $\mathcal R_{\chi_r}(G_r)$.

Recall that $A^{F_r}$ denotes a maximal $F_r$-split torus in $G$, with centralizer $T$, and $B = TU$ is the $F$-rational Borel subgroup defining the dominant Weyl chamber in $X_*(A^{F_r})_{\mathbb R}$ or $X_*(A)_\mathbb R$.  Let $\rho$ denote the half-sum of the $B$-positive roots of $G$.  

Fix any uniformizer $\varpi$ for the field $F$.  Consider a regular dominant cocharacter $\nu \in X_*(A^{F_r})$ and set $u = \nu(\varpi) \in T(F_r)$.  Also, let $\tau = N_r(\nu)$, a regular dominant cocharacter in $X_*(A)$, and set $t = \tau(\varpi) \in T(F)$.  Thus $t = N_r(u)$.  

A construction of Casselman \cite{Cas} and Deligne \cite{Del} associates to any element $g \in G(F)$ an $F$-rational parabolic $P_g = M_gN_g$ with unipotent radical $N_g$ and $F$-rational Levi factor $M_g$.  Following Labesse \cite{Lab90}, we may construct the analogous objects associated to $u\theta \in G_r \rtimes \langle \theta \rangle$.  The resulting subgroups $M_{u\theta}$ (resp. $P_{u\theta}$) of $G$ can be characterized as the set of elements $g \in G$ such that $(u\theta)^n g (u\theta)^{-n}$ remains bounded as $n$ ranges over all integers (resp. all positive integers).  In our situation $M_{u\theta} = M_t = T$, and $P_{u\theta} = P_t = B$.  

Now fix any $F$-Levi subgroup $M$.   Write $G_r$ (resp. $M_r$) for the set of $F_r$-points $G(F_r)$ (resp. $M(F_r)$).  Then $M_r$ acts on $G_r \times M_r$ on the left by $m_1(g,m) = (m_1g, m_1 m \theta(m_1)^{-1})$.  Let $[g,m]$ denote the equivalence class of $(g,m)$ in the quotient space $M_r \backslash G_r \times M_r$.  There is a natural morphism of $p$-adic analytic manifolds
\begin{align} \label{Cas-Del_morph}
M_r \backslash G_r \times M_r & \rightarrow G_r \\
[g,m] & \mapsto g^{-1}m \theta(g). \notag
\end{align}
It is well-known (cf. \cite{Cl90, Lab90}) that this morphism is generically one-to-one and that the normalized absolute value of its Jacobian at the point $[g,m]$ is $|D_{G(F)/M(F)}(\mathcal N(m))|_F$.  

Let us apply this to our Levi subgroup $M_{u\theta} = T$; we are now concerned only with the quotient by $T(F_r)_1$ of the compact open subset $I_r \times T(F_r)_1 u \subset G_r \times T_r$.  The action of $T(F_r)_1$ on $I_r \times T(F_r)_1 u$ is defined by restricting the above-defined action of $T_r$ on $G_r \times T_r$.  
The map
\begin{align} \label{map}
T(F_r)_1 \backslash I_r \times T(F_r)_1 u & \rightarrow G_r \\
[k,mu] &\mapsto k^{-1}mu\theta(k) \notag
\end{align} 
is injective and has the absolute value of its Jacobian everywhere equal to the non-zero number 
\begin{equation} \label{Jac}
|{\rm Jac}_{[k,mu]}|_F = \delta^{-1}_{B_r}(u) = q^{r\langle 2\rho, \nu \rangle}.
\end{equation}  
Thus its image is a compact {\em open} subset, which we will denote $\mathfrak L_u$.

\begin{defn}
We define the elementary function $\phi_{u,\chi_r}$ on $G_r$ to vanish off of $\mathfrak L_u$, and to have value on $\mathfrak L_u$ given by
\begin{equation} \label{elem_def}
\phi_{u,\chi_r}(k^{-1}mu \theta(k)) = \chi_r^{-1}(m)
\end{equation}
for $k \in I_r$ and $m \in T(F_r)_1$.
\end{defn}

We claim that $\phi_{u, \chi_r}$ is well-defined.  Indeed, suppose
$$
k_1^{-1} m_1 u \theta(k_1) = k_2^{-1} m_2 u \theta(k_2),
$$
for $k_1, k_2 \in I_r$ and $m_1, m_2 \in T(F_r)_1$.   By Lemma \ref{Lab_lem}(i) below, the element $k = k_2 k_1^{-1}$ belongs to $T_r$.  Thus it belongs to $T(F_r)_1$ and commutes with $u$, and we have
$$
m_1 k \theta(k^{-1}) = m_2
$$
which implies $\chi_r(m_1) = \chi_r(m_2)$.  Thus we need only the following result due to Labesse \cite{Lab90}.

\begin{lemma} \label{Lab_lem} Suppose $g_1, g_2 \in G(\bar{F})$ and $m_1, m_2 \in T(F_r)_1$.  Extend $\theta$ to an element of ${\rm Gal}(\bar{F}/F)$ and use the same symbol $\theta$ to denote 
the induced automorphism of $G(\bar{F})$.  We have the following statements.
\smallskip

\noindent {\rm (i)}  If $g_1^{-1} m_1 u \theta(g_1) = g_2^{-1} m_2 u \theta(g_2)$, then $g_1 \in T(\bar{F}) g_2$.
\smallskip

\noindent {\rm (ii)}  The $\theta$-centralizer of $m_1u$ is $T$; hence $m_1u$ is strongly 
$\theta$-regular.
\smallskip

\noindent {\rm (iii)}  If $g_1 \in G_r$ and $g_1^{-1} m_1 u \theta(g_1)$ lies in the support of $\phi_{u, \chi_r}$, then $g_1 \in T_r I_r$.  
\end{lemma}

\begin{proof}
Part (i).  Let $g = g_1 g_2^{-1}$. As in Lemma 1 of \cite{Lab90} one proves, using the fact that $u$ commutes with $m_1$ and $m_2$, that $(u\theta)^n g (u\theta)^{-n}$ remains bounded as $n$ ranges over $\mathbb Z$, and thus $g \in M_{u\theta} = T$.  It is enough to let $n$ range only over multiples of $r$.

Let $n$ be any positive integer.  We have
$$
g_1^{-1} (m_1 u \theta)^{rn} g_1 = g_2^{-1}(m_2u\theta)^{rn} g_2,
$$
and thus 
$$
[N_{r}(m_1)N_{r}(u) \theta^{r}]^n = g\, [N_{r}(m_2) N_{r}(u) \theta^{r}]^n\, g^{-1}.
$$
The terms in each set of brackets $[\cdots]$  pairwise commute, so that we deduce
$$
[N_{r}(u) \theta^{r}]^n \, g \, [N_{r}(u) \theta^{r}]^{-n} = N_{r}(m_1)^{-n} \, g \, N_{r}(m_2)^n.
$$
The right hand side remains bounded as $n \in \mathbb Z$.  The left hand side is 
$(u\theta)^{rn} \, g \, (u\theta)^{-rn}$.  This shows $g \in M_{u\theta} = T$ and proves part (i).

Parts (ii) and (iii) follow from part (i). 
\end{proof}

Thus, $\phi_{u, \chi_r}$ is a smooth compactly-supported function on $G_r$, which is invariant under $\theta$-conjugation by elements in $I_r$.   For $r = 1$ (hence $\theta = {\rm id}$ and $u=t$) we get the analogous smooth compactly supported function $f_{t, \chi}$ on $G$, and it is invariant under conjugation by $I$.

For later use in subsection \ref{end_subsec}, we record the following lemma.

\begin{lemma} \label{bi-invar}
The functions $\phi_{u,\chi_r}$ belong to $\mathcal H(G_r,\rho_r)$.
\end{lemma}

\begin{proof}
Labesse shows in \cite{Lab95}, Proposition IV.1.1, that $\mathfrak L_u = I_r u I_r$; thus the support of $\phi_{u, \chi_r}$ is stable under left or right multiplication by $I_r$.  The first step is to check left $I^+_r$-invariance of $\phi_{u, \chi_r}$ (the right invariance is similar).  

Using the normality of $I^+_r$ in $I_r$, it is sufficient to show that an identity of the form
\begin{equation} \label{ident1}
k_1^{-1} i^+ m_1 u \theta(k_1) = k_2^{-1}m_2 u \theta(k_2),
\end{equation}
with $i^+ \in I^+_r$, $k_1,k_2 \in I_r$, and $m_1, m_2 \in T(F_r)_1$, implies that 
\begin{equation} \label{equal}
\chi_r(m_1) = \chi_r(m_2).
\end{equation}
Let $k := k_2 k_1^{-1}$.  Then (\ref{ident1}) gives
\begin{equation} \label{ident2}
m_2^{-1}ki^+ m_1 = u \theta(k) u^{-1}.
\end{equation}  
Now write $k \in I_r$ according to the Iwahori decomposition (\ref{Iwahori_decomp}) as
\begin{equation} \label{ident3}
k = k_- \cdot k_0 \cdot k_+ \in I_{r\bar{U}} \cdot T(F_r)_1 \cdot I_{rU}.
\end{equation}
The element $ki^+$ has an Iwahori decomposition of the form 
\begin{equation} \label{ident4}
ki^+ = k_- i^+_-  \cdot k_0 i^+_0 \cdot k_+i^+_+  
\end{equation}
for suitable elements $i^+_- \in I^+_r \cap \bar{U}$, $i^+_0 \in I^+_r \cap T$ and $i^+_+ \in I^+_r \cap U$.  So equation (\ref{ident2}) has the form
\begin{equation} \label{form}
m_2^{-1}(k_- i^+_-)m_2 \cdot m_2^{-1}(k_0 i^+_0) m_1 \cdot m_1^{-1}(k_+i^+_+) m_1 = u\theta(k_-)u^{-1} \cdot \theta(k_0) \cdot u\theta(k_+) u^{-1}.
\end{equation}
Using the uniqueness of the decompositions in the ``big cell'' $\bar{U} \cdot T \cdot U$, we deduce that
$$
m_2^{-1} (k_0i^+_0)  m_1 = \theta(k_0).
$$
This implies (\ref{equal}) and proves the $I^+_r$-bi-invariance.

To finish the proof, we need to consider an equality of the form
$$
m_0 k_1^{-1} m_1 u \theta(k_1) = k_2^{-1} m_2 u \theta(k_2)
$$
for $k_1,k_2, m_1, m_2$ as above and $m_0 \in T(F_r)_1$, and show that $\chi_r(m_2) = \chi_r(m_0 m_1)$.  We can write this in the form
$$
k_1^{-1} i^+ (m_0m_1 u)\theta(k_1) = k_2^{-1} m_2 u\theta(k_2)
$$
for a suitable element $i^+ \in I^+_r$.  Now the argument above applied to this expression yields $\chi_r(m_0m_1) = \chi_r(m_2)$, as desired.
\end{proof}

The utility of the functions $\phi_{u, \chi_r}$ stems from two principles:

\begin{itemize}
\item The twisted orbital integrals of $\phi_{u, \chi_r}$ are easy to compute; we can thus prove $\phi_{u, \chi_r}$ and $f_{t, \chi}$ are associated functions for all $u$.
\item The character of $\phi_{u, \chi_r}$ on any $\theta$-stable representation $\Pi$ of $G_r$ is easy to compute, and is non-vanishing only on representations $\Pi$ belonging to $\mathcal R_{\chi_r}(G_r)$; thus we will get explicit character identities in the parameter $u$ coming from the matching of $\phi_{u,\chi_r}$ and $f_{t, \chi}$.
\end{itemize}

We will establish these facts in the next few subsections.  The next subsection contains a key lemma about the norm homomorphism.

\subsection{A lemma on the norm homomorphism}

\begin{lemma} \label{norm}  Let $T$ be any $F$-torus.  Then the norm homomorphism $N_r: T(F_r)_1 \rightarrow T(F)_1$ is surjective.
\end{lemma}

This is well-known (and easy) when $T$ is $F$-split; in the general case it can be proved indirectly from the fact that the characteristic functions $1_{T(F_r)_1}$ and $1_{T(F)_1}$ are associated (\cite{Ko86b}).  Our purpose here is only to give a more elementary and direct proof.  In this article we will use this lemma only in the case where $T$ is split over $L$.

\begin{proof}
Following \cite{HR1}, let $\mathcal T^\circ$ denote the neutral component of the lft Neron model $\mathcal T$ associated to $T$.  Then $\mathcal T^\circ$ is a smooth $\mathcal O$-group scheme with connected geometric fibers.  By loc.~cit. we have $T(F)_1 = \mathcal T^\circ(\mathcal O)$ and $T(F_r)_1 = \mathcal T^\circ(\mathcal O_r)$.  

Therefore, we need to prove that the group 
$$
\dfrac{\mathcal T^\circ(\mathcal O_r)^\theta}{N_r(\mathcal T^\circ(\mathcal O_r))} = \widehat{H}^0(F_r/F, \mathcal T^\circ(\mathcal O_r))$$
vanishes.  

First note that for each positive integer $n$, we have
\begin{equation} \label{H1=0}
H^1(F_r/F, \mathcal T^\circ(\mathcal O_r/\varpi^n\mathcal O_r)) = 0.
\end{equation}
Indeed, since inflation is always injective on $H^1$, it is enough to prove the vanishing of $H^1(\langle \theta \rangle, \mathcal T^\circ(\mathcal O_L/\varpi^n\mathcal O_L))$.  But this is  a quotient of the group $H^1(\langle \theta \rangle, \mathcal T^\circ(\mathcal O_L))$ which is well-known to be trivial 
(cf. e.g. \cite{Ko97}, (7.6.1)).

Further, the ${\rm Gal}(F_r/F)$-module $\mathcal T^\circ(\mathcal O_r/\varpi^n\mathcal O_r)$ is finite and thus has trivial Herbrand quotient (cf. \cite{Ser}, VIII, Prop. 8).  Hence from (\ref{H1=0}) we deduce
$$\widehat{H}^0(F_r/F, \mathcal T^\circ(\mathcal O_r/\varpi^n\mathcal O_r)) = 0.$$
Thus to prove the lemma it would be enough to prove that the natural map
\begin{equation} \label{nat_map}
\dfrac{\mathcal T^\circ(\mathcal O_r)^\theta}{N_r(\mathcal T^\circ(\mathcal O_r))} \longrightarrow \underset{n}{\varprojlim} \dfrac{\mathcal T^\circ(\mathcal O_r/\varpi^n\mathcal O_r)^\theta}{N_r(\mathcal T^\circ(\mathcal O_r/\varpi^n\mathcal O_r))}
\end{equation}
is injective.  But this is a straightforward exercise which goes as follows.  Let $t_n$ be the reduction modulo $\varpi^n$ of $t \in \mathcal T^\circ(\mathcal O_r)^\theta$, and suppose for every $n$ we have $t_n = N_r(s_n)$ for some $s_n \in \mathcal T^\circ(\mathcal O_r/\varpi^n\mathcal O_r)$.  Using (\ref{H1=0}) and the surjectivity of the transition maps 
\begin{equation} \label{transitions}
\mathcal T^\circ(\mathcal O_r/\varpi^{n+1}\mathcal O_r) \rightarrow \mathcal T^\circ(\mathcal O_r/\varpi^n \mathcal O_r)
\end{equation}
we may inductively alter the $s_n$'s so that they are compatible under the maps (\ref{transitions}), yet still satisfy $N_r(s_n) = t_n$ for all $n$.  Then $s := \underset{n}{\varprojlim} \, s_n$ is an element in $\mathcal T^\circ(\mathcal O_r)$ with $N_r(s) = t$.  This proves the lemma.
\end{proof}

\subsection{Orbital integrals}
 
Normalize the Haar measure $dg$ on $G_r$ (resp. $ds$ on $T_r$) so that $I_r$ (resp. $T(F_r)_1$) has measure 1.  Use this to define the quotient measures $d\bar{g}$ on $T(F)\backslash G_r$ (resp. $d\bar{s}$ on $T(F)\backslash T_r$).

\begin{prop}\label{orbital_prop}  The following statements hold.  
\smallskip

\noindent {\rm (1)}  For $\delta \in G_r$, ${\rm TO}_{\delta \theta}(\phi_{u, \chi_r})$ is non-zero if and only if $\delta$ is $\theta$-conjugate in $G_r$ to an element of the form $mu$, $\, m \in T(F_r)_1$.  For such $m$ we have
\begin{equation} \label{TO_mu}
{\rm TO}_{mu\theta}(\phi_{u, \chi_r}) = \chi_r^{-1}(m).
\end{equation}
\smallskip

\noindent {\rm (2)}  For $\gamma \in G$, ${\rm O}_\gamma(f_{t,\chi})$ is non-zero if and only if $\gamma$ is conjugate in $G$ to an element of the form $m_0 t$, $\, m_0 \in T(F)_1$.  For such $m_0$ we have
\begin{equation} \label{O_m0}
{\rm O}_{m_0t}(f_{t,\chi}) = \chi^{-1}(m_0).
\end{equation}
\end{prop}

\begin{proof}
It is enough to prove (1).  Let $\delta \in G_r$.  If ${\rm TO}_{\delta \theta}(\phi_{u, \chi_r}) \neq 0$, then $\delta$ is $\theta$-conjugate under $G_r$ to an element of the form $m u$, where $m \in T(F_r)_1$.  Hence we might as well assume $\delta = mu$; in that case $G^\circ_{\delta \theta}(F) = G^\circ_{mu\theta}(F) = T(F)$ (Lemma \ref{Lab_lem}(ii)).  Now if $g \in G_r$ is such that $g^{-1}mu\theta(g) \in {\rm supp}(\phi_{u,\chi_r})$, then by Lemma \ref{Lab_lem}(ii) we have $g \in T_rI_r$.  Then using our choice of measures we get
\begin{align*}
\int_{T(F)\backslash G_r} \phi_{u, \chi_r}(g^{-1}mu\theta(g)) \, d\bar{g} &= \int_{T(F)\backslash T_r} \phi_{u,\chi_r}(s^{-1}mu \theta(s)) \, d\bar{s} \\
&= \int_{T(F)\backslash T_r} \phi^{T_r}_{1,\chi_r}(s^{-1}m\theta(s)) \, d\bar{s}.
\end{align*}
Here $\phi^{T_r}_{1,\chi_r}$ denotes the elementary function on $T_r$ associated to $u = 1$.  This is precisely the characteristic function $1_{T(F_r)_1}$ of $T(F_r)_1$ times the character $\chi_r^{-1}$.  Thus we have proved
\begin{equation} \label{TO_computation}
{\rm TO}^{G_r}_{mu\theta}(\phi_{u, \chi_r}) = \chi^{-1}_r(m) \, {\rm TO}^{T_r}_{m\theta}(1_{T(F_r)_1}). 
\end{equation}
Using 
\begin{equation} \label{second_H1=0}
H^1(F_r/F, T(F_r)_1) = 0
\end{equation}
(see the proof of (\ref{H1=0})), one easily checks that ${\rm TO}^{T_r}_{m\theta}(1_{T(F_r)_1}) = {\rm vol}_{d\bar{s}}(T(F)\backslash T(F) T(F_r)_1)$.  By our choice of measures, this volume equals 1.
\end{proof}

\begin{lemma} \label{slick}  We have the following statements.
\smallskip

\noindent {\rm (i)}  Let $m_0, m'_0 \in T(F)_1$.  Then $m_0 t$ and $m'_0 t$ are stably-conjugate in $G$ if and only if $m_0 = m_0'$.
\smallskip

\noindent {\rm (ii)}  Let $m, m' \in T(F_r)_1$.  Then $m u$ and  $m'u$ are stably $\theta$-conjugate in $G_r$ if and only if they are $\theta$-conjugate by an element of $T(F_r)_1$.
\end{lemma}

\begin{proof}
Part (i) follows from Lemma \ref{Lab_lem}(i), taking $r=1$ (hence $\theta = {\rm id}$ and $u = t$).  Next, $mu$ is stably $\theta$-conjugate to $m'u$ if and only if $N_r(m)t$ is stably-conjugate to $N_r(m')t$ (by \cite{Ko82}, Prop. 5.7).  By part (i), this holds if and only if $N_r(m) = N_r(m')$.  Because of (\ref{second_H1=0}), this holds if and only if $m$ and $m'$ are $\theta$-conjugate in $T(F_r)_1$.
\end{proof}

\begin{cor} \label{O=SO}
For every $\delta \in G_r$, ${\rm SO}_{\delta \theta}(\phi_{u, \chi_r}) = {\rm TO}_{mu \theta}(\phi_{u,\chi_r})$, where $mu$, $\, m \in T(F_r)_1$, is any element of this form which is stably $\theta$-conjugate to $\delta$.
\end{cor}

\begin{proof}
If $\delta$ is not stably-conjugate to an element of the form $mu$, then both sides are zero.  If $\delta = mu$, then the non-zero summands in ${\rm SO}_{\delta \theta}(\phi_{u, \chi_r})$ are of the form ${\rm TO}_{m'u}(\phi_{u,\chi_r})$ (they come with coefficient 1 since $m'u$ is strongly $\theta$-regular).  By Lemma \ref{slick}(ii), only one term appears in this sum.
\end{proof}

\begin{prop} \label{elem_assoc}
The functions $f_{t, \chi}$ and $\phi_{u, \chi_r}$ are associated.
\end{prop}

\begin{proof}
First we assume that ${\rm SO}_{\gamma}(f_{t,\chi}) \neq 0$ and prove this implies that $\gamma$ is a norm from $G_r$.  By Proposition \ref{orbital_prop}(2), we see that $\gamma$ is stably-conjugate in $G$ to an element of the form $m_0t$, $\, m_0 \in T(F)_1$.  By Lemma \ref{norm}, $m_0 = N_r(m)$ for some $m \in T(F_r)_1$.  Since $t = N_r(u)$, this proves that $\gamma$ is a norm from $G_r$.

Next we assume $\gamma$ is the norm of an element $\delta \in G_r$.  If ${\rm SO}_{\delta \theta}(\phi_{u, \chi_r}) \neq 0$, then $\delta$ is stably $\theta$-conjugate to an element of the form $mu$, and then by Corollary \ref{O=SO} and Proposition \ref{orbital_prop} we have
\begin{equation} \label{SO1}
{\rm SO}_{\delta \theta}(\phi_{u, \chi_r}) = {\rm TO}_{mu\theta}(\phi_{r, \chi_r}) = \chi_r^{-1}(m).
\end{equation}
But letting $m_0 := N_r(m)$, \cite{Ko82} Prop. 5.7 shows that $\gamma$ is stably-conjugate to the element $m_0t$ and the same reasoning then gives us
$$
{\rm SO}_\gamma(f_{t,\chi}) = {\rm O}_{m_0t}(f_{t,\chi}) = \chi^{-1}(m_0),
$$
which coincides with (\ref{SO1}).

If ${\rm SO}_{\delta \theta}(\phi_{u, \chi_r}) = 0$, then $\delta$ is not stably $\theta$-conjugate to an element of the form $mu$.  But then $\gamma$ is not stably conjugate to an element of the form $m_0 t$, and thus ${\rm SO}_\gamma(f_{t,\chi}) = 0$ as well.  Indeed, if $\gamma$ were stably conjugate to $m_0t$, then writing $m_0 = N_r(m)$ for some $m$ (Lemma \ref{norm}), we would have that $N_r(\delta)$ is stably conjugate to $N_r(mu)$.  But then by \cite{Ko82}, Prop. 5.7, it would be the case that $\delta$ is stably $\theta$-conjugate to $m u$, contrary to what we saw above.
\end{proof}

\subsection{Traces of elementary functions}

Let $\Pi$ denote a $\theta$-stable admissible representation of $G(F_r)$, and fix an intertwiner 
$I_\theta : \Pi ~ \widetilde{\rightarrow} ~ \Pi^\theta$.  Let $\Theta_{\Pi \theta}$ denote the locally integrable function of Harish-Chandra representing the functional $\phi \mapsto \langle {\rm trace} \, \Pi I_\theta, \phi\rangle $, for $\phi \in C^\infty_c(G(F_r))$.  Thus,
\begin{equation} \label{character_defn}
\langle {\rm trace} \, \Pi I_\theta, \phi\rangle = \int_{G_r} \Theta_{\Pi \theta}(g) \phi(g) \, dg.
\end{equation}

Following \cite{Lab90},$\S3$, we calculate this trace for $\phi = \phi_{u,\chi_r}$ using the map (\ref{map}).  Let $d\bar{k}$ denote the quotient measure on $T(F_r)_1\backslash I_r$ denoted $d\bar{g}$ earlier.   The change of variable formula yields
\begin{align} 
\langle {\rm trace}\,  \Pi I_\theta, \phi_{u,\chi_r} \rangle &= \int_{T(F_r)_1\backslash I_r} 
\int_{T(F_r)_1} |{\rm Jac}_{[k,mu]}|_F \cdot \Theta_{\Pi \theta}(k^{-1} mu\theta(k))\cdot \phi_{u, \chi_r}(k^{-1}mu\theta(k)) \, dm \, d\bar{k} \notag \\
&= q^{r\langle 2\rho, \nu \rangle} \int_{T(F_r)_1} \Theta_{\Pi \theta}(mu) \, \chi_r^{-1}(m) \, dm \notag \\
&= q^{r\langle 2\rho, \nu \rangle} \int_{T(F_r)_1} \Theta_{\Pi_U\theta}(mu) \, \chi_r^{-1}(m) \, dm. \label{ch}
\end{align}
Here $\Pi_U$ is the Jacquet module of $\Pi$ corresponding to the Borel subgroup $B_r  = T_r U_r$, and the equality $\Theta_{\Pi \theta}(mu) = \Theta_{\Pi_U \theta}(mu)$ we used is the twisted version due to Rogawski \cite{Rog} of a theorem of Casselman \cite{Cas}.

\begin{lemma} \label{trace_cond}
Suppose $\Pi$ is an irreducible $\theta$-stable object in $\mathcal R(G_r)$.  If $\langle {\rm trace}\,  \Pi I_\theta, \phi_{u,\chi_r} \rangle \neq 0$, then $\Pi$ belongs to the subcategory $\mathcal R_{\chi_r}(G_r)$.  
\end{lemma}

\begin{proof}
The non-vanishing of (\ref{ch}) shows that $\Pi_U^{\chi_r} \neq 0$.  Then by (\ref{Jac_isom}), we see $\Pi^{\rho_r} \neq 0$.  Now the result follows by Proposition \ref{type}.
\end{proof} 

\subsection{Computation of traces of elementary functions} \label{trace_calc}

In this subsection, we write $W$ (resp. $W_r$) for the relative Weyl group associated to the $F$-split (resp. $F_r$-split) torus $A$ (resp. $A^{F_r}$) in $G$.   

For any irreducible $\theta$-stable representation $\Pi$ of $G_r$ with intertwiner $I_\theta$, we will compute here the trace
$$
\langle {\rm trace} \, \Pi I_\theta, \phi_{u,\chi_r} \rangle.
$$  
The computation is similar\footnote{In the final displayed equation on p.526 of \cite{Lab90}, the calculation of ${\rm trace}(I_{\tilde{\lambda}})_{N_0}(u\theta)$ is done assuming a particularly simple form for the intertwiner on the induced module $I_{\tilde \lambda}$ (more precisely, on its Jacquet module $(I_{\tilde \lambda})_{N_0})$).  Unless $\tilde \lambda$ is regular, there will be other possible intertwiners on $(I_{\tilde \lambda})_{N_0}$.   The formula in Proposition 7' of loc.~cit. implicitly assumes that the intertwiner on $\Pi$ is inherited from such a simple intertwiner on a suitable induced module $I_{\tilde \lambda}$.  The statement of our Proposition \ref{trace_prop} is more complicated because we do not make any assumptions about the intertwiner on $\Pi$.} to that in \cite{Lab90}.   By Lemma \ref{trace_cond}, we may assume $\Pi$ belongs to the category $\mathcal R_{\chi_r}(G_r)$ (otherwise the trace is zero).  Suppose the supercuspidal support of $\Pi$ is $(T_r, \xi')_{G_r}$ for some extension $\xi'$ of a $W_r$-conjugate of $\chi_r$.  

Let $\Xi$ denote the set of characters on $T(F_r)$ which extend some $W_r$-conjugate of $\chi_r$.   Thus, $\Xi$ contains $W_r\xi'$.  Let $\Xi(\chi_r) \subset \Xi$ consist of those whose restriction to $T(F_r)_1$ is precisely $\chi_r$.  Let $\Xi^\theta$ resp. $\Xi(\chi_r)^\theta$ denote the subset of $\theta$-fixed elements in $\Xi$ resp. $\Xi(\chi_r)$.  It will become clear below that the trace above is zero unless $\Xi^\theta(\chi_r) \neq \emptyset$.  From now on, $\xi'$ will denote an arbitrary element of $\Xi$, not just the one we started with.

Now we fix further notation related to $\Pi$.  Recall that $\Pi_U$ denotes the Jacquet module of $\Pi$ relative to $U_r$.  Because $(\Pi^\theta)_U = (\Pi_U)^\theta$, the intertwiner $I_\theta: \Pi ~ \widetilde{\rightarrow} ~ \Pi^\theta$ induces an intertwiner $I_\theta: \Pi_U ~ \widetilde{\rightarrow} ~ \Pi^\theta_U$.  Recall that $\Pi_U$ is a subquotient of 
$$
(i^{G_r}_{B_r}(\xi'))_U = \delta_{B_r}^{1/2} \bigoplus_{w \in W_r} \mathbb C_{^w\xi'}
$$
(cf. \cite{Cas1}), where $\mathbb C_{^w\xi'}$ is the 1-dimensional representation of $T_r$ corresponding to the character $^w\xi'$.  Thus there is a well-defined subset $\Xi(\Pi) \subset \Xi$ and positive multiplicities $a_{\xi',\Pi} =: a_{\xi'}$ for $\xi' \in \Xi(\Pi)$ such that, as $T_r$-representations,
\begin{equation} \label{Pi_U_form}
\Pi_U = \delta_{B_r}^{1/2} \bigoplus_{\xi' \in \Xi(\Pi)} \mathbb C_{\xi'}^{a_{\xi'}}.
\end{equation}
Since $\Pi^\theta_U$ is isomorphic to $\Pi_U$, the set $\Xi(\Pi)$ is stabilized by $\theta$.  


For the following statement and proof, we set  $\Xi(\Pi)^\theta = \Xi^\theta \cap \Xi(\Pi)$, and $\Xi(\Pi,\chi_r)^\theta := \Xi(\chi_r)^\theta \cap \Xi(\Pi)$.  Further, 
for $\xi' \in \Xi^\theta$, we set 
$${\rm tr}(I_\theta, \Pi,\xi') := 
\langle {\rm trace} \,I_\theta \, ; \, \delta_{B_r}^{1/2} \mathbb C^{a_{\xi'}}_{\xi'} \rangle.$$  Finally, for $\xi' \in \Xi(\chi_r)$, we may write 
$$\xi' = \widetilde{\chi}^\varpi_r \, \eta'$$ 
for a unique unramified character $\eta'$ on $T_r$.

\begin{prop} \label{trace_prop}
In the notation above, we have  
\begin{equation} \label{trace_eq}
\langle {\rm trace} \, \Pi I_\theta, \phi_{u, \chi_r} \rangle = q^{\langle \rho, \tau \rangle} 
\sum_{\xi' \in \Xi(\Pi,\chi_r)^\theta}\, \eta'(u) \, {\rm tr}(I_\theta, \Pi,\xi').
\end{equation}
\end{prop}

\begin{proof} 
We will rewrite the integral in (\ref{ch}).  Using (\ref{Pi_U_form}) it takes the form
$$
\int_{T(F_r)_1} \chi_r^{-1}(m) \, \langle {\rm trace} \, \Pi_U (mu)I_\theta \, ; \, \delta_{B_r}^{1/2} \bigoplus_{\xi' \in \Xi(\Pi)} \mathbb C_{\xi'}^{a_{\xi'}} \rangle \, dm.
$$
The isotypical component $\delta^{1/2}_{B_r} \mathbb C_{\xi'}^{a_{\xi'}}$ can contribute to the trace of $\Pi_U I_\theta$ only if $\xi' \in \Xi(\Pi)^\theta$, and in that case $\Pi_U I_\theta$ preserves that component.  
Thus the integral can be expressed as
\begin{equation} \label{form1}
\sum_{\xi' \in \Xi(\Pi)^\theta} \int_{T(F_r)_1} \chi_r^{-1}(m) \, \langle {\rm trace} \, \Pi_U(mu)I_\theta \, ;\, \delta_{B_r}^{1/2} \mathbb C_{\xi'}^{a_{\xi'}} \rangle \, dm.
\end{equation} 
On $\delta_{B_r}^{1/2}\mathbb C_{\xi'}^{a_{\xi'}}$ appearing here, $\Pi_U(mu)$ acts by the scalar
$$
\delta_{B_r}^{1/2}(u) \, \xi'(mu) = q^{-\langle \rho, \tau \rangle} \, \xi'(m) \, \eta'(u).
$$
(We used that $\widetilde{\chi}^\varpi_r(u) =1 $.)  Thus (\ref{form1}) becomes
\begin{equation*}
q^{-\langle \rho, \tau \rangle} \sum_{\xi' \in \Xi(\Pi)^\theta} \, \eta'(u) \, {\rm tr}(I_\theta,\Pi,\xi')\,  \int_{T(F_r)_1} \chi^{-1}_r(m) \, \xi'(m) \, dm.
\end{equation*}
The integral on the right hand side is non-vanishing (and equal to 1) if and only if $\xi'$ equals $\chi_r$ on $T(F_r)_1$, that is, if and only if $\xi' \in \Xi(\Pi,\chi_r)^\theta$.  The result now follows.  
 \end{proof}

We shall also need the following result.  Recall the idempotent $e_{\rho} \in \mathcal 
Z(G,\rho)$ defined to have support $I$ and to take value $\rho(k)^{-1}$ at $k \in I$.

\begin{lemma}
In the notation above, we have
\begin{equation} \label{tr_e_rho_calc}
\sum_{\xi' \in \Xi(\Pi,\chi_r)^\theta} {\rm tr}(I_\theta,\Pi, \xi') = \langle {\rm trace}\, \Pi I_\theta \, , \, e_{\rho_r} \rangle.
\end{equation}
\end{lemma}

\begin{proof}
The right hand side can be written as 
$$
\int_{I_r} \langle {\rm trace}\, \Pi(g)I_\theta \, , \, \Pi^{\rho_r} \rangle\,  e_{\rho_r}(g) \, dg = \langle {\rm trace}\, I_\theta \,, \, \Pi^{\rho_r} \rangle.$$  
By  (\ref{Jac_isom}) and (\ref{Pi_U_form}), this can be written as
$$
\langle {\rm trace}\, I_\theta \,,\, \delta_{B_r}^{1/2} \bigoplus_{\xi' \in \Xi(\Pi,\chi_r)^\theta} \mathbb C_{\xi'}^{a_{\xi'}} \rangle,
$$
which is the left hand side in (\ref{tr_e_rho_calc}).  
\end{proof}

\begin{cor} \label{trace_u=1_cor}
We have
\begin{equation} \label{trace_u=1_eq}
\langle {\rm trace} \, \Pi I_\theta, \phi_{u,\chi_r} \rangle|_{u=1} = \langle {\rm trace} \, \Pi I_\theta, e_{\rho_r} \rangle.
\end{equation}
In particular, if $r=1$, we have
\begin{equation} \label{trace_u=1_r=1}
\sum_{\xi \in \Xi(\pi,\chi)} \, {\rm dim}(\mathbb C_\xi^{a_{\xi,\pi}}) = \langle {\rm trace} \, \pi, e_\rho \rangle.
\end{equation}
\end{cor}

\section{Proof in the strongly regular elliptic case} \label{ell_proof}

\subsection{Local data adapted to $\mathcal R_\chi(G)$} \label{local_data_sec}

In this subsection, we assume $G$ is any unramified group over $F$ with $G = G_{\rm ad}$.  Let ${\rm Irr}_\chi(G)$ (resp. ${\rm Irr}^\theta_{\chi_r}(G_r)$) denote the set of irreducible (resp. irreducible $\theta$-stable) admissible representations in $\mathcal R_\chi(G)$ (resp. $\mathcal R_{\chi_r}(G_r)$).

Inspired by Hales \cite{Ha}, we define {\em local data adapted to $\mathcal R_\chi(G)$} to consist of the data (a),(b), and (c), subject to conditions (1) and (2) below:

(a) An indexing set $\mathcal I$ (possibly infinite);

(b) A collection of complex numbers $a_i(\pi)$ for $i \in \mathcal I$ and $\pi \in {\rm Irr}_\chi(G)$;

(c) A collection of complex numbers $b_i(\Pi)$ for $i \in \mathcal I$ and $\Pi \in {\rm Irr}^\theta_{\chi_r}(G_r)$.

\smallskip

(1) For $i$ fixed, the constants $a_i(\pi)$ and $b_i(\Pi)$ are zero for all but finitely many $\pi$ and $\Pi$.

(2) For $\phi \in \mathcal H(G_r,\rho_r)$ and $f \in \mathcal H(G,\rho)$ the following are equivalent:

\hspace{.5in}  (A)  For all $i$, we have $\sum_{\pi} a_i(\pi) \langle {\rm trace} \, \pi, f  \rangle = \sum_{\Pi} b_i(\Pi) \langle {\rm trace} \, \Pi I_\theta , \phi \rangle$;

\hspace{.5in}  (B)  For all strongly regular elliptic semisimple norms $\gamma = \mathcal N(\delta)$, we have $${\rm SO}_{\gamma}(f) = {\rm SO}_{\delta \theta}(\phi).$$

\smallskip

One can prove that such local data exist using a result of Clozel (cf. \cite{H09a}, Proposition 8.3.1), in exactly the same way as for parahoric Hecke algebras (see loc.~cit.~$\S8.5$).  

A word about the intertwiners $I_\theta$: these are local intertwiners coming from canonically defined intertwiners at the global (adelic) level (cf. \cite{H09a}, Proposition 8.3.1).  But the global-to-local process involves some choices and hence the $I_\theta$ here are not canonically defined.  However, since $\Pi$ is irreducible, $I_\theta$ is defined up to scalar, and such scalars can be absorbed into the coefficients $b(\Pi)$.  Hence we are free to normalize the $I_\theta$ however we like. 

\subsection{End of proof} \label{end_subsec}

Assume $G = G_{\rm ad}$.  Fix $\phi \in \mathcal Z(G_r,\rho_r)$. and $f = b_r(\phi) \in \mathcal Z(G,\rho)$.  

We now complete the proof of Theorem \ref{main_thm}.  By Lemma \ref{wlog_norm}, we may assume $\gamma$ is a norm.  By our reduction section, we may assume $\gamma$ is a strongly regular elliptic semisimple norm, say $\gamma = \mathcal N(\delta)$.  To check the identity ${\rm SO}_{\delta \theta}(\phi) = {\rm SO}_\gamma(f)$, we follow the method of Labesse \cite{Lab90}: we use the associated elementary functions $\phi_{u,\chi_r}, f_{t,\chi}$ to give, as $u$ ranges, sufficiently many character identities to establish (2)(A) in the local data for the pair $\phi, f$. (Recall $u = \varpi^\nu$ and $t = N_r(u) = \varpi^{\tau}$, and that $\nu$ ranges over the semigroup of regular dominant cocharacters 
$X_*(A^{F_r})^{++}$ in the group $X_*(A^{F_r})$).

We need to show that for each $i$, the identity in (2)(A) above holds.  Fixing $i$ and dropping it from our notation, we need to prove that
\begin{equation} \label{char_id}
\sum_{\pi} a(\pi) \, \langle {\rm trace}~\pi , b_r(\phi) \rangle =  \sum_{\Pi} b(\Pi) \, \langle {\rm trace}~\Pi I_\theta , \phi \rangle.
\end{equation}
By Lemma \ref{bi-invar} and Proposition \ref{elem_assoc}, the equivalent statements (2)(A) and (2)(B) hold for each pair $\phi_u := \phi_{u,\chi_r}$, $f_t := f_{t,\chi}$.  Therefore
\begin{equation} \label{elem_id}
\sum_{\pi} a(\pi) \, \langle {\rm trace}~\pi , f_t \rangle = \sum_{\Pi} b(\Pi) \, \langle {\rm trace}~\Pi I_\theta , \phi_u \rangle.
\end{equation} 

Let us rewrite the left hand side using (\ref{trace_eq}) in the case where $r=1$.  We get
\begin{equation} \label{lhs}
\sum_{\pi} a(\pi) \, \langle {\rm trace}~\pi , f_t \rangle ~ = ~ q^{\langle \rho, \tau \rangle} \, \sum_{\xi_0 \in \Xi(\chi)/W_\chi} \, \sum_{\pi \in i^G_B(\xi_0)} \, a(\pi) \, \sum_{\xi \in \Xi(\pi,\chi)} \eta(t) \, ({\rm dim}\, \mathbb C_\xi^{a_{\xi,\pi}}).
\end{equation}
Here $\xi_0$ ranges over any set of representatives for the $W_\chi$-orbits on the set $\Xi(\chi)$ of characters on $T(F)$ extending $\chi$.   Since $\eta(t) = \eta_r(u)$, (\ref{lhs}) is a linear combination of characters which are norms evaluated on $u = \varpi^\nu$.

We would like to write down a similar expression for the right hand side of (\ref{elem_id}) and then compare it to (\ref{lhs}).  An apparent difficulty is that given two different elements in $\Xi(\chi_r)^\theta$ which are $W_{\chi_r}$ conjugate, it could happen that one of them is a norm and the other is not.  To remedy this, we will group the elements according to $W_\chi$-conjugacy rather than $W_{\chi_r}$-conjugacy (note that 
$W_\chi \subseteq W^\theta_{\chi_r}$).  The following lemma shows that this avoids the problem just mentioned.  We leave the proof to the reader.

\begin{lemma} \label{W_chi_conj}
If $\xi'_1, \xi'_2$ belong to $\Xi(\chi_r)^\theta$ and are $W_\chi$-conjugate, then one is a norm if 
and only if the other one is.
\end{lemma}

\medskip
Now again using (\ref{trace_eq}), we can write the right hand side of (\ref{elem_id}) in the following way:
\begin{align} \label{rhs}
\sum_{\Pi} b(\Pi) \, \langle {\rm trace}  \, \Pi I_\theta \, , \, & \phi_u \rangle  ~~ =  ~ ~  \notag \\
q^{\langle \rho, \tau \rangle} \, \sum_{\xi'_0 \in \Xi(\chi_r)^\theta/W_{\chi}} \,& \sum_{\overset{\Pi \, \mbox{\tiny s.t.}}{W_\chi\xi'_0 \cap \Xi(\Pi,\chi_r)^\theta \neq \emptyset}} \,  b(\Pi) \, \sum_{\overset{\xi' \in}{W_\chi \xi'_0 \cap \Xi(\Pi,\chi_r)^\theta}} \eta'(u) \, {\rm tr}(I_\theta, \Pi, \xi'). 
\end{align}
Let $(\ref{rhs})_n$ (resp. $(\ref{rhs})_{nn}$) denote the contribution to (\ref{rhs}) coming from those $\xi'_0$ which are norms (resp. not norms).  By Lemma \ref{W_chi_conj}, all (resp. none) of the terms $\eta'$ appearing in $(\ref{rhs})_n$ (resp. $(\ref{rhs})_{nn}$) are norms.  

Now we regard each side of the identity 
\begin{equation}\label{(=)}
(\ref{lhs}) = (\ref{rhs})_n + (\ref{rhs})_{nn}
\end{equation}
 as a linear combination of characters on $\nu \in X_*(A^{F_r})^{++}$.  The independence of characters asserts that any non-empty set of distinct characters on $X_*(A^{F_r})$ is linearly independent over $\mathbb C$.  The proof shows that the set of their restrictions to $X_*(A^{F_r})^{++}$ remains independent.  The equation (\ref{(=)}) holds for all regular dominant $\nu$, thus by this linear independence statement, it must hold for $\nu = 0$, that is, for $u = 1$.  In the self-evident notation, this means that we have two identities
\begin{align} \label{2identities}
(\ref{lhs})|_{u=1} &= (\ref{rhs})_{n}|_{u=1} \\
0 &= (\ref{rhs})_{nn}|_{u=1}.
\end{align}
In fact the linear independence gives something stronger: the contributions of each $\xi'_0 \in \Xi(\chi_r)^\theta/W_\chi$, to the above equations satisfy corresponding identities.  Making use of Corollary \ref{trace_u=1_cor}, this may be stated as follows.

\begin{lemma} \label{sep_lem_1}
For each $\xi'_0 \in \Xi(\chi_r)^\theta/W_\chi$, we have  
\begin{align} \label{sep_eq_1}
\sum_{\overset{\xi_0 \in \Xi(\chi)/W_\chi}{\mbox{\tiny s.t.} \, \xi_{0r} = \xi'_0}} \, \sum_{\pi \in i^G_B(\xi_0)} a(\pi) \, \langle {\rm trace} \, \pi \, , \, & e_\rho \rangle ~~ = ~~ \\
\sum_{\overset{\Pi \, \mbox{\tiny s.t.}}{W_\chi \xi'_0 \cap \Xi(\Pi,\chi_r)^\theta \neq \emptyset}} \, &b(\Pi) \, \sum_{\overset{\xi' \in}{W_\chi\xi'_0 \cap \Xi(\Pi,\chi_r)^\theta}}  \, {\rm tr}(I_\theta, \Pi, \xi'). \notag
\end{align}
\end{lemma} 
Note that if $\xi'_0$ is not a norm, the left hand side of (\ref{sep_eq_1}) vanishes, and the right hand does as well by our linear independence argument.  

Next we fix any $\xi'_0 \in \Xi(\chi_r)^\theta$ and sum this equation over all $W_\chi$-orbits in $W_{\chi_r}\xi'_0 \cap \Xi(\chi_r)^\theta$.  In that sum, the right hand side will simplify, because of the identity (\ref{tr_e_rho_calc}).  Indeed, if $W_{\chi_r}\xi'_0$ contains a norm, we may assume $\xi'_0$ itself is a norm, and the $W_\chi$-orbits in $W_{\chi_r}\xi'_0$ which are not norms contribute nothing to the right hand side of (\ref{sep_eq_1}), and hence we get an equality
\begin{align} \label{CaseA}
\sum_{\overset{\xi_0 \in \Xi(\chi)/W_\chi}{\mbox{\tiny s.t.}\, \xi_{0r} \in W_{\chi_r}\xi'_0}} \, \sum_{\pi \in i^G_B(\xi_0)} \, & a(\pi) \, \langle {\rm trace}\, \pi \, , \, e_\rho \rangle \\
~ = ~ \sum_{\Pi \in i^{G_r}_{B_r}(\xi'_0)} \, & b(\Pi) \, \langle {\rm trace} \, \Pi I_\theta \, , \, e_{\rho_r} \rangle. \notag
\end{align}
On the other hand, if $W_{\chi_r}\xi'_0$ contains no norms, then (\ref{CaseA}) remains true: the left hand side obviously vanishes, and the right hand side vanishes as a result of the linear independence discussion above.  

Recall that if $\xi'_0 = \xi_{0r}$, then $b(\phi)$ acts on $i^G_B(\xi_0)^{\rho}$ by the same scalar by which $\phi$ acts on $i^{G_r}_{B_r}(\xi'_0)^{\rho_r}$ (Lemma \ref{left_ch_br_lem}).  Thus multiplying (\ref{CaseA}) by $ch_{\xi_0}(b(\phi)) = ch_{\xi_{0r}}(\phi)$ if $\xi'_0 = \xi_{0r}$, and by $ch_{\xi'_0}(\phi)$ in the case where $\xi'_0$ is not a norm, we get (\ref{CaseA}) but with $e_{\rho_r}$ resp. $e_\rho$ replaced by $\phi$ resp. $b(\phi)$.  Then summing the resulting formula over all $\xi'_0 \in \Xi(\chi_r)/W_{\chi_r}$ yields the desired formula (\ref{char_id}).  This completes the proof of Theorem \ref{main_thm}. \qed

\begin{Remark}
The results of this article are closely related to at least four earlier works on the fundamental lemma, namely \cite{Cl90}, \cite{Lab90}, \cite{Ha}, \cite{H09a}.  We point out that all of these rely crucially on a result of Keys \cite{Keys}.  It is interesting to note that the present article does not need Keys' result at all.  The reason is that we are using Labesse's method of elementary functions: in \cite{Lab90}, Keys' result is used only to separate, in the character identities appearing in the local data, the spherical representation from other representations in the same unramified principal series representation.  There is no need for any such separation argument in this article.
\end{Remark} 

\section{Base change for $I^+$-level} \label{bc_I+_sec}

In this section we will construct a base change homomorphism
\begin{equation} \label{bc_I+}
b_r: \mathcal Z(G_r,I^+_r) \rightarrow \mathcal Z(G,I^+)
\end{equation}
and show that the fundamental lemma also holds for it, as claimed in Corollary \ref{bc_I+_cor}. 


First we abuse notation by writing $T(k_r)$ in place of $\mathcal T^\circ(k_r)$, where $\mathcal T^\circ$ is the neutral component of the lft N\'{e}ron model $\mathcal T$ for the $F$-torus $T$.  Thus there is a canonical isomorphism
$$
T(k_r) = T(F_r)_1/T(F_r)^+_1
$$
and the set of depth-zero characters on $T(F_r)_1$ can be canonically identified with $T(k_r)^\vee$, the set of characters of the finite group $T(k_r)$.  For $\chi' \in T(k_r)^\vee$, we consider the corresponding Hecke algebra $\mathcal H(G_r, I_r, \chi')$ and its center $\mathcal Z(G_r,I_r, \chi')$.

Recall the idempotents $e_{\chi'} \in \mathcal Z(G_r,I_r, \chi')$: the function $e_{\chi'}$ is supported on $I_r$ and takes value ${\rm vol}(I_r)^{-1} \rho_{\chi'}(x)^{-1}$ at $x \in I_r$.  These form a complete set of orthogonal idempotents in $\mathcal H(G_r,I^+_r)$ (cf.~e.g.~\cite{HR2}, proof of Prop. 12.1.1).  Hence there is an algebra {\em monomorphism}
\begin{align} \label{mono}
\mathcal Z(G_r, I^+_r) &\hookrightarrow \prod_{\chi' \in T(k_r)^\vee} \mathcal Z(G_r, I_r, \chi') \\
z &\mapsto (ze_{\chi'})_{\chi'}. \notag
\end{align}
We need to identify the image of (\ref{mono}).  As it clearly suffices, we discuss this in the case $r=1$ and make self-evident adjustments to our notation.
\begin{lemma} \label{image_lem}
The image of (\ref{mono}) consists of the elements $(z_{\chi})_{\chi}$ with the following property: for any relative Weyl group element $w \in W$, any $\chi \in T(k)^\vee$, and any extension $\tilde{\chi}: T(F) \rightarrow \mathbb C^\times$ of $\chi$, the scalar by which $z_{\chi}$ acts on $i^{G}_{B}(\widetilde{\chi})^{\rho_{\chi}}$ coincides with the scalar by which $z_{^w\chi}$ acts on $i^{G}_{B}(\widetilde{\chi})^{\rho_{\,^w\chi}}$. 
\end{lemma}

\begin{proof}
It is clear that elements in the image of (\ref{mono}) satisfy the property in question.  We assume $(z_\chi)_\chi$ satisfies the property and show that it comes from a $z \in \mathcal Z(G,I^+)$.

Fix $\chi_0$ and consider its $W$-orbit $[\chi_0]$.  Recall (\cite{HR2}, $\S12$) the {\em central} idempotent $e_{[\chi_0]} := \sum_{\chi \in [\chi_0]} e_\chi$ in $\mathcal H(G,I^+)$.  This is the projector onto the Bernstein component $\mathcal R_{\chi_0}(G)$, in the sense that for every smooth representation $V$ of $G(F)$, the $G(F)$-module spanned by $e_{[\chi_0]} V$ is the component of $V$ belonging to $\mathcal R_{\chi_0}(G)$.  Using Proposition \ref{type}, we deduce
\begin{equation} \label{this_module}
\bigoplus_{\chi \in [\chi_0]} i^G_B(\widetilde{\chi_0})^{\rho_\chi} = i^G_B(\widetilde{\chi_0})^{I^+} = e_{[\chi_0]} i^G_B(\widetilde{\chi}_0),
\end{equation}
and the $G$-module generated by the latter is just $i^G_B(\widetilde{\chi_0})$.  By our hypothesis, the element $\sum_{\chi \in [\chi_0]} z_\chi$ operates by a scalar $ch_{\widetilde{\chi_0}}$ on (\ref{this_module}).  As $\widetilde{\chi_0}$ varies these scalars give a regular function on the Bernstein variety $\mathfrak X_{\chi_0}$.  Let $\mathcal Z(G)$ denote the Bernstein center of $G(F)$, viewed as 
certain distributions on $C^\infty_{c}(G(F))$ or equally as regular functions on the union of all Bernstein varieties.  We see there is an element $z_{[\chi_0]} \in \mathcal Z(G) * e_{[\chi_0]} \subset \mathcal Z(G,I^+)$ such that $z_{[\chi_0]} e_{\chi} = z_\chi$ for each $\chi \in [\chi_0]$.  Now let $z := \sum_{\chi_0 \in W\backslash T(k)^\vee} z_{[\chi_0]} \in \mathcal Z(G,I^+)$.  Then $z$ goes to $(z_\chi)_\chi$ under (\ref{mono}).
\end{proof}

If $\chi'$ is a norm, then $\chi' = \chi_r := \chi \circ N_r$ for some $\chi \in T(k)^\vee$, and since $N_r$ is surjective, the character $\chi$ is uniquely determined by $\chi'$.  Now consider the composition of the projection onto the factors indexed by norms
$$
\prod_{\chi' \in T(k_r)^\vee} \mathcal Z(G_r, I_r, \chi') \rightarrow \prod_{\chi \in T(k)^\vee} \mathcal Z(G_r, I_r, \chi_r)
$$
followed by the map
$$
\prod_{\chi \in T(k)^\vee} \mathcal Z(G_r,I_r, \chi_r) \rightarrow \prod_{\chi \in T(k)^\vee} 
\mathcal Z(G,I,\chi)
$$
which is given in each factor by the homomorphism $b_r$ of Definition \ref{bc_chi_defn}.

\begin{lemma}
This composition $\prod_{\chi'} \mathcal Z(G_r,I_r, \chi') \rightarrow \prod_{\chi} \mathcal Z(G,I,\chi)$ above takes the image of (\ref{mono}) into the corresponding image for $r=1$.
\end{lemma}

\begin{proof}
Let $z \in \mathcal Z(G_r,I^+_r)$.  We need to check that the tuple $(b_r(ze_{\chi_r}))_{\chi}$ satisfies the property of Lemma \ref{image_lem}.  By Lemma \ref{left_ch_br_lem}, the scalar by which $b_r(ze_{\chi_r})$ acts on $i^G_B(\chi)^{\rho_\chi}$ is the scalar by which $ze_{\chi_r}$ acts on $i^{G_r}_{B_r}(\chi_r)^{\rho_{\chi_r}}$.  Apply this to both $\chi$ and $^w\chi$, for $w$ any element of $W$, the relative Weyl group over $F$.
\end{proof}

This lemma allows us to define the base change homomorphism for $I^+$-level.

\begin{defn} \label{I+_bc_def}
The homomorphism $b_r: \mathcal Z(G_r,I^+_r) \rightarrow \mathcal Z(G,I^+)$ is the unique homomorphism making the following diagram commute:
$$
\xymatrix{
\mathcal Z(G_r,I^+_r) \ar[r] \ar[d]_{b_r} & \prod_{\chi'} \mathcal Z(G_r,I_r, \chi') \ar[d] \\
\mathcal Z(G,I^+) \ar[r] & \prod_{\chi} \mathcal Z(G,I,\chi).}
$$
\end{defn}

The following lemma explains why it is permissible to forget all the factors corresponding to $\chi'$ which are not norms.

\begin{lemma} \label{chi'_TO_van}
Suppose $\chi' \in T(k_r)^\vee$ is not a norm.  Then all twisted orbital integrals at $\theta$-semisimple elements vanish on $\mathcal H(G_r, I_r, \chi')$.
\end{lemma}

\begin{proof}
Fix $\phi \in \mathcal H(G_r, I_r, \chi')$ and consider any $i \in I_r$.  The substitution $g \mapsto gi$ in ${\rm TO}_{\delta \theta}(\phi) = \int_{G^\circ_{\delta\theta}\backslash G_r} \phi(g^{-1}\delta \theta(g)) \, d\bar{g}$ yields
$$
{\rm TO}_{\delta \theta}(\phi) = (\chi' \theta(\chi')^{-1})(i) \, {\rm TO}_{\delta \theta}(\phi).
$$
Thus if the twisted orbital integral does not vanish, we have $\chi' = \theta(\chi')$.  But then $\chi'$ is a norm: the $\theta$-invariant elements in $T(k_r)^\vee$ are precisely the norms.
\end{proof}

With these remarks in hand, it is easy to see how Corollary \ref{bc_I+_cor} follows from Theorem \ref{main_thm}. \qed

\section{Corrigenda to \cite{H09a}}  \label{corr_sec}

All notation will be that of \cite{H09a}, but note the correction in notation discussed below.  

\subsection{\cite{H09a}, section 2.2}  In section 2.2 of \cite{H09a}, there is a minor misstatement.  It has no effect on the main results of \cite{H09a}, but nevertheless this corrigendum seems necessary in order to avoid potential confusion.  I am very grateful to Brian Smithling and Tasho Kaletha, who brought this misstatement to my attention.

\subsubsection{Correction}

Here is the precise misstatement.  In \cite{H09a}, section 2.2, the ``ambient'' group scheme $\mathcal G_{{\bf a}_J}$ was incorrectly identified with the group scheme whose group of $\mathcal O_L$-points is the full fixer of the facet ${\bf a}_J$.   In the notation of Bruhat-Tits \cite{BT2}, which I intended to follow in \cite{H09a}, the group scheme whose group of $\mathcal O_L$-points is the full fixer of ${\bf a}_J$ is denoted 
$\widehat{\mathcal G}_{{\bf a}_J}$.  The group scheme $\widehat{\mathcal G}_{{\bf a}_J}$ is defined and characterized in this way in \cite{BT2}, 4.6.26-28.  

The group scheme denoted $\mathcal G_{{\bf a}_J}$ is defined in loc.~cit.~4.6.26 (cf. also 4.6.3-6).  In general, it can be a bit smaller than $\widehat{\mathcal G}_{{\bf a}_J}$ (see below).  

The correction is the following: in \cite{H09a}, the symbol $\mathcal G_{{\bf a}_J}$ should now be interpreted as the group denoted by this symbol in \cite{BT2}, a potentially proper subgroup of the full fixer $\widehat{\mathcal G}_{{\bf a}_J}$.  

We still have, as stated in \cite{H09a}, (2.3.2) and (2.3.3), the equalities
\begin{align} 
J(L) = \mathcal G^\circ_{{\bf a}_J}(\mathcal O_L) &= T(L)_1 \cdot \mathfrak U_{{\bf a}_J}(\mathcal O_L) \label{2.3.2} \\
\mathcal G_{{\bf a}_J}(\mathcal O_L) &= T(L)_b \cdot \mathfrak U_{{\bf a}_J}(\mathcal O_L). \label{2.3.3}
\end{align}
In general,
$$
\mathcal G^\circ_{{\bf a}_J}(\mathcal O_L) = \widehat{\mathcal G}^\circ_{{\bf a}_J}(\mathcal O_L) \subset \mathcal G_{{\bf a}_J}(\mathcal O_L) \subset \widehat{\mathcal G}_{{\bf a}_J}(\mathcal O_L),
$$
and both inclusions can be strict.  

\subsubsection{Clarifications}  We discuss the effect of the above correction on subsequent statements in \cite{H09a}.

\smallskip

${\bf 1.}$  Theorem 2.3.1 of \cite{H09a} remains valid as stated, but can be slightly augmented: equation (2.3.1) can be replaced by
\begin{equation} \label{thm_2.3.1}
J(L) = {\rm Fix}({\bf a}_J^{\rm ss}) \cap G(L)_1 = \mathcal G_{{\bf a}_J}(\mathcal O_L) \cap G(L)_1 = \widehat{\mathcal G}_{{\bf a}_J}(\mathcal O_L) \cap G(L)_1.
\end{equation}
Cf. \cite{HR1}, Remark 11.

\smallskip

${\bf 2.}$  Contrary to \cite{H09b}, line above equation (2.3.2), our ${\mathcal G}_{{\bf a}_J}$ should not now be identified with the scheme $\widehat{\mathcal G}_{{\bf a}^{\rm ss}_J}$ of \cite{BT2}.

\smallskip

${\bf 3.}$  Corollary 2.3.2 of \cite{H09b} remains valid, with the same proof. Indeed, when $G_L$ is split we have $T(L)_b = T(\mathcal O_L) = T(L)_1$ and then from (\ref{2.3.2}) and (\ref{2.3.3}) above we see that $\mathcal G^\circ_{{\bf a}_J}(\mathcal O_L) = \mathcal G_{{\bf a}_J}(\mathcal O_L)$. 

\smallskip

${\bf 4.}$  Lemma 2.9.1 of \cite{H09b} remains valid as stated, but in the proof (especially in equations (2.9.1) and (2.9.2)) the symbols $\mathcal G_{{\bf a}_J}(\mathcal O_L)$ and $\mathcal G_{{\bf a}^M_J}(\mathcal O_L)$ should be replaced by  $\widehat{\mathcal G}_{{\bf a}_J}(\mathcal O_L)$ and $\widehat{\mathcal G}_{{\bf a}^M_J}(\mathcal O_L)$, respectively.

\subsubsection{Example}

It is sometimes but usually not the case that $\mathcal G_{{\bf a}_J}(\mathcal O_L) = \widehat{\mathcal G}_{{\bf a}_J}(\mathcal O_L)$.  The following is perhaps the simplest example where this equality fails\footnote{Brian Smithling and Tasho Kaletha provided me with another example for the split group ${\rm SO}(2n)$.}.  Take $G$ to be the split group ${\rm PSp}(4)$, and let ${\bf a}_J$ denote the non-special vertex in a base alcove.  Then let $\tau$ denote the element in the stabilizer $\Omega \subset \widetilde{W}(L)$ of the base alcove, which interchanges the two special vertices and fixes ${\bf a}_J$.  The element $\tau$ does not belong to the group $\mathcal G^\circ_{{\bf a}_J}(\mathcal O_L) = \mathcal G_{{\bf a}_J}(\mathcal O_L)$ (cf.~${\bf 3}$ above), since $\tau$ does not belong to $G(L)_1$.  On the other hand $\tau \in \widehat{\mathcal G}_{{\bf a}_J}(\mathcal O_L)$ since it fixes ${\bf a}_J$ and $G(L)^1 = G(L)$ (cf. \cite{BT2}, 4.6.28).

\subsection{\cite{H09a}, Lemma 2.9.1(b)}  The proof of the equality $J \cap P = (J \cap M)(J \cap N)$ given in Lemma 2.9.1(b) is flawed.  Instead, one can deduce this from the following equality
$$
J = (J \cap B)(J \cap \overline{N})(J \cap N)(J \cap M).
$$
which is a consequence of \cite{BT2}, 5.2.4.

\small
\bigskip
\obeylines
\noindent
University of Maryland
Department of Mathematics
College Park, MD 20742-4015 U.S.A.
email: tjh@math.umd.edu

\end{document}